\documentclass{article}%
\usepackage{amsfonts}
\usepackage{amsmath}
\usepackage{amssymb}
\usepackage{graphicx}%
\newtheorem{theorem}{Theorem}
\newtheorem{acknowledgement}[theorem]{Acknowledgement}

\newtheorem{condition}[theorem]{Condition}

\newtheorem{corollary}[theorem]{Corollary}
\newtheorem{criterion}[theorem]{Criterion}
\newtheorem{definition}[theorem]{Definition}

\newtheorem{lemma}[theorem]{Lemma}

\newtheorem{proposition}[theorem]{Proposition}
\newtheorem{remark}[theorem]{Remark}

\begin{document}

\title{\thanks{This work is partially supported by the Institute of Mathematical
Sciences, CUHK}Shafarevich's Conjecture for CY Manifolds I \\(Moduli of CY Manifolds)}
\author{Kefeng Liu\\UCLA, Department of Mathematics\\Los Angeles, CA 90095\\and Center of Mathematical Sciences\\Zheijiang University
\and Andrey Todorov\\UC, Department of Mathematics\\Santa Cruz, CA 95064\\Institute of Mathematics\\Bulgarian Academy of Sciences
\and Shing-Tung Yau\\Harvard University\\Department of Mathematics\\Cambridge, Mass. 02138
\and Kang Zuo\\Chinese University of Hong Kong\\Department of Mathematics\\Shatin, Hong Kong}
\maketitle

\begin{abstract}
In this paper we first study the moduli spaces related to Calabi-Yau
manifolds. We then apply the results to the following problem. Let $C$ be a
fixed Riemann surface with fixed finite number of points on it. Given a CY
manifold with fixed topological type, we consider the set of all families of
CY manifolds of the fixed topological type over $C$ with degenerate fibres
over the fixed points up to isomorphism. This set is called Shafarevich set.
The analogue of Shafarevich conjecture for CY manifolds is for which
topological types of CY the Shafarevich set is finite. It is well-known that
the analogue of Shafarevich conjecture is closely related to the study of the
moduli space of polarized CY manifolds and the moduli space of the maps of
fixed Riemann surface to the coarse moduli space of the CY manifolds. We prove
the existence of the Teichm\"{u}ller space of CY manifolds together with a
universal family of marked CY manifolds. From this result we derive the
existence of a finite cover of the coarse moduli space which is a non-singular
quasi-projective manifold. Over this cover we construct a family of polarized
CY manifolds. We study the moduli space of maps of the fixed Riemann with
fixed points on it to the moduli space of CY manifolds constructed in the
paper such that the maps map the fixed points on the Riemann surface to the
discriminant locus. If this moduli space of maps is finite then Shafarevich
conjecture holds. We relate the analogue of Shafarevich problem to the
non-vanishing of the Yukawa coupling. We give also a counter example of the
Shafarevich problem for a class of CY manifolds.

\end{abstract}
\listoftables
\tableofcontents

\section{Introduction}

We define the Teichm\"{u}ller space of any complex manifold as the quotient of
all integrable complex structures on M by the action of the group of
diffeomorphisms isotopic to the identity. One of the main achievements in the
theory of Riemann surfaces was the construction of the Teichm\"{u}ller spaces
of Riemann surfaces of different genus. This was done by Teichm\"{u}ller, L.
Bers, L. Alfors. The definition of a Teichm\"{u}ller space implies that the
mapping class group of the Riemann surface of a fixed genus acts on it. The
quotient is the moduli space. In the case of Riemann surfaces of genus greater
than or equal to one, it is a well known fact that the Teichm\"{u}ller space
is a domain of holomorphy. Important role in the study of moduli of Riemann
surfaces is played by the Weil-Petersson metric. In case of Riemann surfaces
of genus greater than one it was noticed by H. Mazur that the Weil-Petersson
metric is not complete. It was Mumford who proved that the moduli space of
Riemann surfaces is a quasi-projective variety. He used geometric invariant
theory. Later the same result was proved by Bers by using the Weil-Petersson metric.

In case the Riemann surface is an elliptic curve, i.e. of genus 1, then it is
a well known fact that the Teichm\"{u}ller space is the upper half plane and
the Weil-Petersson metric in this case is the Poincare metric on it.

Very little is known about the Teichm\"{u}ller space of higher dimensional
K\"{a}hler manifolds. In the case of polarized abelian varieties of complex
dimension $g\geq2,$ the Teichm\"{u}ller space is the Siegel upper half plane
of genus $g.$ In case of polarized algebraic K3 surfaces, the Teichm\"{u}ller
space is an open and everywhere dense subset in
\[
\mathfrak{h}_{2,19}:=SO_{0}(2,19)/SO(2)\times SO(19).
\]
The complement to the Teichm\"{u}ller space in\ $\mathfrak{h}_{2,19}$
corresponds to algebraic K3 surfaces admitting double rational points. We will
show that the Teichm\"{u}ller space $\widetilde{T}$(M) of a Calabi-Yau
manifold M exists and it has finite number of irreducible components. Each
component $\mathcal{T}$(M) of $\widetilde{T}$(M) is a non-singular complex
manifold. Moreover over $\widetilde{T}$(M) there exists a universal family of
marked polarized CY manifolds up to the action of a finite group of complex
analytic automorphisms which acts trivially on the middle cohomology.

The existence of the coarse moduli space $\mathfrak{M}_{L}$(M) of polarized CY
manifolds as a quasi-projective variety was proved by Viehweg. See \cite{V}.
Yau conjectured that Mumford-Chow stability of the canonically embedded
projective manifolds is equivalent to the existence of K\"{a}hler Eimstein
metric. Recently Donaldson proved that constant scalar curvature implies
stability. Combining the result of Donaldson with the solution of the Calabi
conjecture given by Yau, one can reprove the result of Viehweg.

In this paper we will prove the existence of a finite cover $\mathcal{M}_{L}%
$(M) of the coarse moduli space $\mathfrak{M}_{L}$(M) of polarized CY
manifolds such that $\mathcal{M}_{L}$(M) is a non-singular quasi-projective
variety and over $\mathcal{M}_{L}$(M) there exists a family of polarized CY
manifolds. The construction of $\mathcal{M}_{L}$(M) and the family over it is
based on the existence of family of CY manifolds over the Hilbert scheme. One
of the main difficulty in using the method of Hilbert schemes is the existence
of the group $G$ of complex analytic automorphisms of a CY manifold which acts
trivially on the middle cohomology. It turned out that $G$ acts on all CY
which are fibres over the connected component of the Hilbert scheme that
contains as a fibre the initial CY manifold. We also examine some relations of
the existence of the group $G$ with the failure of the global Torelli Theorem
for CY manifolds.

We will apply the results of this paper to the analogue of the Shafarevich's
conjecture about finiteness of the families of polarized CY manifolds over a
fixed Riemann surface with prescribed points of degenerate fibres. In his
article published in the Proceedings of the International Congress of
Mathematicians, Stockholm meeting held in 1962, Shafarevich wrote:

\textit{''One of the main theorems on algebraic numbers connected with the
concept of discriminant is Hermit's theorem, which states that the number of
extensions }$k^{\prime}/k$ \textit{of a given degree and given discriminant is
finite. This theorem may be formulated as follows: the number of extensions
}$k^{\prime}/k$\textit{\ of a given degree whose critical prime divisors
belong to a given finite set S is finite.''}

Inspired by this result of Hermit, Shafarevich conjectured in \cite{sh}:
\textit{\textquotedblright There exists only a finite number of fields of
algebraic functions }$K/k$\textit{\ of a given genus }$g\geq1$\textit{, the
critical prime divisors of which belong to a given finite set }$S."$

In one unpublished work, Shafarevich proved his conjecture in the setting of
hyperelliptic curves. On page 755 of \cite{sh}, in his remarks on his papers,
he wrote: \textit{\textquotedblright Here two statements made in lecture are
mixed into one: formulation of a result and a conjecture. The result was
restricted to the case of hyperelliptic curves while the conjecture concerned
general curves...This conjecture became much more attractive after A. N.
Parshin proved that it implies the Mordell conjecture...In 1983 it was proved
by G. Faltings (Invent. Math. \textbf{73, }}439-466(1983))\textit{, with a
proof of the Mordell conjecture as a consequence.\textquotedblright}

The formulation of the Shafarevich problem is the following. Let C be a fixed,
non-singular algebraic curve, and let E be a fixed effective divisor on C such
that all the points in E have multiplicity 1. Define Sh(C,E,Z) to be the set
of all isomorphism classes of projective varieties
\[
\pi:\mathcal{Z\rightarrow}C
\]
with a fibre of \textquotedblright type" $Z$ such that the only singular
fibres are over the set E. The general Shafarevich type problem
is:\textbf{\ \textquotedblright}\textit{For which \textquotedblright
type\textquotedblright\ of varieties }$Z$ \textit{and data (C,E)\ is such that
}Sh (C,E,Z) \textit{is finite} ?

Previous works on Shafarevich type problems include the following results. In
the case when $Z$ is a curve of genus $g>1$ and E is empty, Parshin proved in
\cite{Par} that Sh (C,E,Z) is finite, and Arakelov proved in \cite{Ar} the
finiteness in the case E is not empty. Faltings constructed examples showing
that Sh (C,E,Z) is infinite for abelian varieties of dimension $\geq8.$ See
\cite{Falt}. Saito an Zucker extended the construction of Faltings to the
setting when Z is an algebraic polarized K3 surface. They were able to
classify all cases when the set Sh (C,E,Z) is infinite. They were not
considering polarized families. See \cite{SZ}. Faltings proved the
Shafarevich's conjecture over number fields and thus he proved the Mordell
conjecture. Yves Andre proved the analogue of Shafarevich's conjectures over
the number fields for K3 surfaces. See \cite{an}. Using techniques from
harmonic maps Jost and Yau analyzed Sh (C,E,Z) for a large class of varieties.
See \cite{JY}. Ch. Peters studied finiteness theorems by considering
variations of Hodge structures and utilizing differential geometric aspects of
the period map and the associated metrics on the period domain. See \cite{Pe}.

A. Parshin and E. Bedulev have proved the boundedness for families of
algebraic surfaces over a fixed algebraic curve assuming that all the fibres
are non-singular. E. Bedulev and E. Viehweg have proved the boundedness for
families of algebraic surfaces over a fixed algebraic curve, which possibly
have singular fibres. See \cite{BW}.

Migliorni, Kov\'{a}c and Zhang proved that any family of minimal algebraic
surfaces of general type over a curve of genus g and m singular points such
that 2g-2+m$\leq0$ is isotrivial. See \cite{Ko}, \cite{Mi}, \cite{Z} and
\cite{BW}.

Recently very important results of Viehweg and the last named author appeared
in \cite{VZ}, \cite{VZ1} and \cite{VZ2}. Brody hyperbolicity was proved for
the moduli space of canonically polarized complex manifolds. In \cite{VZ2} 6.2
a) They proved the boundedness for Sh(C,E,Z) for arbitrary $Z,$ with
$\omega_{Z}$ semi-ample. They also established that the automorphism group of
moduli stacks of polarized manifolds is finite. The rigidity property for the
generic family of polarized manifolds was also proved. The basis idea in the
proof is to consider the non-vanishing property of the maximal iterated
Kodaira-Spencer map. One notices that in the case of Calabi-Yau manifolds this
maximal iterated Kodaira-Spencer map is related to the Yukawa coupling. In
this paper the relation between the Yukawa coupling and rigidity problems for
CY manifolds is first formally formulated.

In this paper we are going to study the analogue of the Shafarevich's
conjecture about finiteness of the families of polarized CY manifolds over a
fixed Riemann surface with prescribed points of degenerate fibres. Any family
of polarized CY\ manifolds induces a map from the base without discriminant
locus to the coarse moduli space. Thus Shafarevich conjecture is reduced to
prove that moduli space of maps of Riemann surface into the moduli space of
polarized CY manifolds with some additional properties is a finite set. First
we noticed that one can replace the coarse moduli space with a finite cover
$\mathcal{M}_{L}$(M). Using the local Torelli Theorem and a result of
Griffiths and Schmid, the finite cover $\mathcal{M}_{L}$(M) of the moduli
space of polarized CY manifolds admits such a metric. Once the existence of a
K\"{a}hler metric with non-negative curvature is constructed \cite{Lu}, it is
not difficult to show that the rigidity of any family of CY manifolds over a
Riemann surface with fixed points of degenerations, Yau's form of Schwarz
Lemma and the Bishop compactness imply Shafarevich conjecture. We proved the
analogue of the Shafarevich's conjecture for CY manifolds by using the Yukawa
coupling to check rigidity:

\begin{condition}
\label{LTYZ}Suppose that
\begin{equation}
\pi:\ \mathcal{X}\rightarrow\text{C} \label{FAM}%
\end{equation}
is a family of polarized CY manifolds over a Riemann surface C. Let let
$t_{0}\in$C such that for M=$\pi^{-1}(t_{0})$ the following condition is
satisfied; For any non zero $\phi\in H^{1}($M,$T^{1,0}),$ $\wedge^{n}\phi
\neq0$ in $H^{n}($M,$\left(  \Omega_{\text{M}}^{n}\right)  ^{\ast})$. Then the
family $\left(  \ref{FAM}\right)  $ is rigid.
\end{condition}

The condition \ref{LTYZ} brings up an interesting relation between rigidity
and mirror symmetry. We would also like to compare Condition \ref{LTYZ} with
some of the results obtained in \cite{VZ1}. In \cite{VZ3} some important
applications of the Condition \ref{LTYZ} are obtained. In \cite{VZ3} the
finiteness of those families, whose iterated Kodaira-Spencer maps have the
same length as the length of the iterated Kodaira-Spencer of the moduli space,
are obtained. E. Viehweg and the last named author constructed in \cite{VZ4}
rigid families of CY hypersurfaces for which the condition \ref{LTYZ} does not
hold, i.e. the Yukawa coupling is zero. They also have constructed non-rigid
families of hypersurfaces in $\mathbb{CP}^{n}$ for any $n>2$ and of any degree
$d>3.$ See \cite{VZ4}. Recently Y. Zhang obtained some important results more
precisely he proved that any Lefshetz family of CY manifold is rigid. See
\cite{Zh} and \cite{Zh1}.

In this paper our approach is close to that of Jost and Yau. Instead of
harmonic maps in this paper we are using holomorphic ones. See \cite{JY}.

\begin{acknowledgement}
This work was finished at the Institute of Mathematical Sciences of the
Chinese University of Hong Kong. We are grateful to the Institute of
Mathematical Sciences of the Chinese University of Hong Kong for their support
and encouragement. Special thanks to Yi Zhang for a careful reading of an
earlier version of our paper. The first and the second author would like to
thank Center of Mathematical Sciences of Zhejiang University where the final
version of this paper was written.
\end{acknowledgement}

\section{Moduli of Polarized CY\ Manifolds}

\subsection{Automorphisms of CY Manifolds that Act Trivially on $H^{n}%
($M,$\mathbb{Z})$}

\begin{theorem}
\label{auto0}Let%
\[
\pi:\ \mathcal{X}\rightarrow\mathcal{K}%
\]
be the Kuranishi family of a polarized CY manifold M=$\pi^{-1}(\tau_{0}),$
$\tau_{0}\in\mathcal{K}$. Suppose that G is a group of holomorphic
automorphisms of M such that G acts trivially on $H_{n}($M,$\mathbb{Z})$ and
preserves the polarization$.$ Then G is a finite group of holomorphic
automorphisms of all the fibres of the Kuranishi family.
\end{theorem}

\textbf{Proof: }Since G acts trivially on $H^{n}($N,$\mathbb{Z})$ and fixes
the polarization class $L$ then the uniqueness the Calabi-Yau metric that
corresponds to $L$ implies that the group $G$ is a group of isometries of the
CY metric$.$ Since on CY manifolds there does not exist global holomorphic
vector fields, we can conclude that $G$ is a discrete subgroup of the
orthogonal group. The compactness of the orthogonal group implies that $G$ is
finite. The local Torelli Theorem implies that
\[
\mathcal{K\subset}\mathbb{P}(H^{n}(\text{M,}\mathbb{Z})\otimes\mathbb{C}).
\]
So G acts on $\mathcal{K}$ and fixes the point $\tau_{0}.$ Since G acts
trivially on $H^{n}($M,$\mathbb{Z})$, then it will act trivially on
$\mathcal{K}$. Next we are going to prove that if $g\in G$ is any element of
$G$ then it acts as complex analytic automorphism on each M$_{\tau}=\pi
^{-1}(\tau)$ for any $\tau\in\mathcal{K}.$ This means that
\begin{equation}
g^{\ast}(I_{\tau})=I_{\tau}, \label{KSK1}%
\end{equation}
where
\[
I_{\tau}\in C^{\infty}(\text{M},Hom(T_{\text{M}}^{\ast},T_{\text{M}}^{\ast
})),\text{ }I_{\tau}^{2}=-id
\]
is the integrable complex structure operator that defines M$_{\tau}$ for
$\tau\in\mathcal{K}$. $T_{\text{N}}^{\ast}$ is the cotangent $C^{\infty}$ real
bundle. It is a well known fact that
\begin{equation}
I_{\tau}=\left(  A_{\tau}\right)  ^{-1}\circ I_{0}\circ A_{\tau}, \label{KSK2}%
\end{equation}
where
\[
A_{\tau}=\left(
\begin{array}
[c]{cc}%
id & \phi(\tau)\\
\overline{\phi(\tau)} & id
\end{array}
\right)  .
\]
See \cite{To89}. Let us recall that if $\Theta_{\text{M}}$ is the holomorphic
tangent bundle then
\[
\phi(\tau)\in C^{\infty}\left(  \text{M},Hom\left(  \Omega_{\text{M}}%
^{1,0},\overline{\Omega_{\text{M}}^{1,0}}\right)  \right)  \approxeq
C^{\infty}\left(  \text{M},\Theta_{\text{M}}\otimes\Omega_{\text{M}}%
^{0,1}\right)
\]
and $\phi(\tau)$ satisfies the equation that guarantees the integrability of
the complex structure operator $I_{\tau}$ defined by $\left(  \ref{KSK2}%
\right)  $:%
\begin{equation}
\overline{\partial}\phi(\tau)=\frac{1}{2}[\phi(\tau),\phi(\tau)]\text{ and
}\overline{\partial}^{\ast}\phi(\tau)=0. \label{KSK3}%
\end{equation}
For all the details see \cite{To89}. Here $\overline{\partial}^{\ast}$ means
the conjugate of the operator $\overline{\partial}$ with respect to the
Calabi-Yau metric corresponding to the polarization class $L.$ If we prove
that for each $g\in G$%
\begin{equation}
g^{\ast}(\phi(\tau))=\phi(\tau) \label{KSK4}%
\end{equation}
then $\left(  \ref{KSK4}\right)  $ implies $\left(  \ref{KSK1}\right)  $ and
so $\phi(\tau)$ is a complex analytic automorphism of N$_{\tau}.$

\textbf{Proof of }$\left(  \ref{KSK4}\right)  $\textbf{:} According to
\cite{To89} if we fix a basis $\phi_{1},...,\phi_{N}$ of harmonic forms of
$\mathbb{H}^{1}($M,$\Theta_{\text{M}})$ then the solution of the equations
$\left(  \ref{KSK3}\right)  $ are given by the power series%
\[
\phi(\tau)=%
{\displaystyle\sum\limits_{i=1}^{N}}
\phi_{i}\tau^{i}+%
{\displaystyle\sum\limits_{i_{1}+...+i_{N}=m>1}}
\phi_{i_{1},...,i_{N}}\left(  \tau^{1}\right)  ^{i_{1}}\times...\times\left(
\tau^{N}\right)  ^{i_{N}}=%
{\displaystyle\sum\limits_{m=1}^{N}}
\phi_{\lbrack m]}(\tau),
\]
where%
\begin{equation}
\phi_{\lbrack m]}(\tau)=%
{\displaystyle\sum\limits_{i_{1}+...+i_{N}=m}}
\phi_{i_{1},...,i_{N}}\left(  \tau^{1}\right)  ^{i_{1}}\times...\times\left(
\tau^{N}\right)  ^{i_{N}} \label{KSK5}%
\end{equation}
and $\phi(\tau)$ satisfies the recurrent relation:%
\begin{equation}
\phi(\tau)=%
{\displaystyle\sum\limits_{i=1}^{N}}
\phi_{i}\tau^{i}+\frac{1}{2}\left(  \overline{\partial}\right)  ^{\ast}%
\circ\mathbb{G}[\phi(\tau),\phi(\tau)]. \label{KSK5a}%
\end{equation}
See \cite{To89}. Here $\mathbb{G}$ is the Green operator associated with the
Laplacian with respect to the CY metric associated with the polarization class
$L.$ Notice that since $G$ is the group of isometries of CY metric then for
any $g\in G$ the Green operator will be invariant, i.e.%
\[
g^{\ast}\mathbb{G=G}.
\]
In \cite{To89} it is proved that if $\phi\in\mathbb{H}^{1}($N,$\Theta
_{\text{N}})$ is a harmonic form with respect to the CY metric then
$\phi\lrcorner\omega_{\text{M}}$ will be a harmonic form of type $(n-1,1).$
This fact together with the fact that $G$ acts trivially $H^{n}($%
M,$\mathbb{C})$ imply that the group acts trivially on $\mathbb{H}^{1}%
($M,$\Theta_{\text{M}}).$ This implies that the linear term of $\left(
\ref{KSK5}\right)  $ satisfies
\begin{equation}
g^{\ast}\left(  \phi_{\lbrack1]}(\tau)\right)  =\phi_{\lbrack1]}(\tau).
\label{KSK6}%
\end{equation}
The proof of the fact that
\[
g^{\ast}\phi(\tau)=\phi(\tau)
\]
is done by induction on the homogeneity of the terms of the power series
$\left(  \ref{KSK5}\right)  .$ Formula $\left(  \ref{KSK6}\right)  $ shows
that $\phi_{\lbrack1]}(\tau)$ is invariant under the action of $g$. Since for
the higher order terms of $\left(  \ref{KSK5}\right)  $ the relation $\left(
\ref{KSK5a}\right)  $ implies%
\[
\phi_{\lbrack k]}=%
{\displaystyle\sum\limits_{p+q=k}}
[\phi_{\lbrack p]}(\tau),\phi_{\lbrack q]}(\tau)]
\]
and $g^{\ast}(\phi_{\lbrack p]})=\phi_{\lbrack p]}$ for $p<k$ then we can
conclude that $g^{\ast}\phi_{\lbrack k]}(\tau)=\phi_{\lbrack k]}(\tau)$ for
any $k>0.$ So the relation$\left(  \ref{KSK4}\right)  $ is proved$.$ Theorem
\ref{auto0} is proved. $\blacksquare$

\subsection{Construction of the Teichm\"{u}ller Space of Marked Polarized CY
Manifolds}

\begin{definition}
\label{Teich}We will define the Teichm\"{u}ller space $\mathcal{T}$(M) of a
compact complex manifold M as follows:
\[
\mathcal{T}(\text{M}):=\mathcal{I}(\text{M})/Diff_{0}(\text{M}),
\]
\textit{where}\
\[
\mathcal{I}(\text{M}):=\left\{  \text{all integrable complex structures on
M}\right\}
\]
\textit{and } Diff$_{0}$(M) \textit{is the group of diffeomorphisms isotopic
to identity. The action of the group Diff(M}$_{0})$ \textit{is defined as
follows; Let }$\phi\in$Diff$_{0}$(M) \textit{then }$\phi$ \textit{acts on
integrable complex structures on M by pull back, i.e. if }
\[
I\in C^{\infty}(\text{M},Hom(T(\text{M}),T(\text{M})),
\]
\textit{then we define } $\phi(I_{\tau})=\phi^{\ast}(I_{\tau}).$
\end{definition}

\begin{definition}
\label{mark}We will call a pair (M; $\gamma_{1},...,\gamma_{b_{n}}$) a marked
CY manifold if M is a CY manifold and $\{\gamma_{1},...,\gamma_{b_{n}}\}$ is a
basis of $H_{n}$(M,$\mathbb{Z}$)/Tor.
\end{definition}

\begin{remark}
Let $\mathcal{K}$ be the Kuranishi space. It is easy to see that if we choose
a basis of $H_{n}$(M,$\mathbb{Z}$)/Tor in one of the fibres of the Kuranishi
family
\[
\pi:\ \mathcal{M\rightarrow K}%
\]
then all the fibres will be marked, since as a $C^{\infty}$ manifold
$\mathcal{X}_{\mathcal{K}}\approxeq$M$\times\mathcal{K}$.
\end{remark}

Next we are going to construct a universal family of polarized marked CY
manifolds
\begin{equation}
\pi:\ \mathcal{U}_{L}\mathcal{\rightarrow T}_{L}(\text{M}) \label{UF}%
\end{equation}
up to the action of the group of complex analytic automorphisms $G$ which acts
trivially on $H^{2}($M,$\mathbb{Z})$ on the fibres$.$ The construction of the
family $\left(  \ref{UF}\right)  $ of marked polarized CY manifolds follows
the ideas of Piatetski-Shapiro and I. R. Shafarevich. See \cite{PS}.

\begin{theorem}
\label{teich}There exists a family of marked polarized CY manifolds
\begin{equation}
\pi:\mathcal{U}_{L}\mathcal{\rightarrow T}_{L}(\text{M}), \label{fam2}%
\end{equation}
which possesses the following properties: \textbf{A. }$\mathcal{T}_{L}($M$)$
is a smooth manifolds of complex dimension $h^{n-1,1}.$ \textbf{B. }The
holomorphic tangent space $\Theta_{\tau,\mathcal{T}_{L}(\text{M})}$ at each
point $\tau\in\mathcal{T}_{L}($M$)$ is naturally identified with $H^{1}%
($M$_{\tau},\Theta_{\text{M}_{\tau}})$ and \textbf{C. }Let
\[
\pi_{\mathcal{C}}:\ \mathcal{Y\rightarrow C}%
\]
be any complex analytic family of marked polarized CY manifolds such that
there exists a point $x_{0}\in\mathcal{C}$ and the fibre $\left(
\pi_{\mathcal{C}}\right)  ^{-1}(x_{0})=$M$_{x_{0}}$ as a marked polarized CY
manifold is isomorphic to some fibre of the family $\left(  \ref{fam2}\right)
$. Then there exists a unique holomorphic map of families
\begin{equation}
\kappa:\ (\mathcal{Y\rightarrow C})\rightarrow\left(  \mathcal{U}%
_{L}\rightarrow\mathcal{T}_{L}(\text{M})\right)  \label{Kod}%
\end{equation}
defined up to a biholomorphic map $\phi$ of the fibres which induces the
identity map on $H_{n}($M$,\mathbb{Z}).$ The restriction of the map $\kappa$
on the base $\mathcal{C}$ is unique.
\end{theorem}

\textbf{Proof: }In this paragraph we will use the following result of Y.-T.
Siu. See \cite{Dem} and \cite{Siu}.

\begin{theorem}
\label{Siu}Let X be an algebraic variety and let $\mathcal{L}$ be an ample
line bundle then $\left(  \mathcal{K}_{\text{X}}\right)  ^{\otimes2}%
\otimes\left(  \mathcal{L}\right)  ^{\otimes m}$ is very ample for any
\begin{equation}
m\geq2+\binom{3n+1}{n}, \label{siu1}%
\end{equation}
where $n$ is the dimension of X.
\end{theorem}

Based on Theorem \ref{Siu} we will prove the following Theorem:

\begin{theorem}
\label{Siu1}Suppose that M is a fixed projective manifold. Let us fix its
cohomology ring $H^{\ast}($M,$\mathbb{Z})$ over $\mathbb{Z},$ its Chern
classed and the polarization class $L\in H^{2}($M,$\mathbb{Z})$ which is the
Chern class of an ample line bundle. Then there are a finite number of
components of the Hilbert scheme that parametrizes all polarized projective
manifolds with fixed cohomology ring $H^{\ast}($M,$\mathbb{Z}),$ Chern classes
and the polarization class $L.$
\end{theorem}

\textbf{Proof: }Let us fix the following data; a projective manifold M with a
canonical class zero, its cohomology ring $H^{\ast}($M,$\mathbb{Z})$ over
$\mathbb{Z},$ its Chern classed and the polarization class $L\in H^{2}%
($M,$\mathbb{Z})$ which is the Chern class of an ample line bundle. According
to a Theorem of Sullivan there are only finite number of $C^{\infty}$
structures on M if the real dimension of M is greater or equal to 5 with the
data mentioned above. See \cite{Sul}. Suppose that $m$ satisfy the inequality
$\left(  \ref{siu1}\right)  $ then from the Riemann-Roch-Hirzebruch Theorem,
the fact that for very ample line bundles $\mathcal{L}$ we have $H^{k}%
($X,$\mathcal{L}$)=0 for $k>0$ and Theorem \ref{Siu}, we can deduce that if we
consider a CY manifold with fixed Chern classes $c_{2},...,c_{n},$ fixed
cohomology ring $H^{\ast}($M,$\mathbb{Z})$ and polarization class $L,$ then
all these K\"{a}hler manifolds with a canonical class zero can be embedded in
a fixed projective space $\mathbb{CP}^{k}.$ Indeed the Hirzebruch-Riemann-Roch
theorem implies that for all K\"{a}hler manifolds with canonical class zero
with fixed cohomology ring $H^{\ast}($M,$\mathbb{Z})$, fixed Chern classes and
the line bundle $\mathcal{L}$ with a fixed ample Chern class $L$ the Euler
characteristics$:$%
\[
\chi(\text{M,}\mathcal{L}^{m})=\dim_{\mathbb{C}}H^{0}(\text{M,}\mathcal{L}%
^{m})=%
{\displaystyle\int\limits_{\text{M}}}
Td(\text{M)}Ch(\mathcal{L}^{m})
\]
have one and the same Hilbert polynomial $\chi($M,$\mathcal{L}^{m})$ for all
such CY manifolds. Using the theory of Hilbert schemes of Grothendieck
developed in \cite{SGA}, we can conclude that there are a finite number of
components of the Hilbert scheme that parametrizes all polarized K\"{a}hler
manifolds with fixed cohomology ring $H^{\ast}($M,$\mathbb{Z}),$ Chern classes
and the polarization class $L.$ $\blacksquare$

In \cite{JT96} we proved the following Theorem:

\begin{theorem}
\label{jt}Let $\mathcal{H}_{L}$ be the Hilbert scheme of \ non-singular CY
manifolds embedded by the linear system $|L^{m}|$ defined by the polarization
class $L,$ then $\mathcal{H}_{L}$ is a non-singular quasi-projective variety.
\end{theorem}

We know from the results in \cite{SGA} that there exists a family of polarized
CY manifolds:
\begin{equation}
\mathcal{Y}_{L}\rightarrow\mathcal{H}_{L}, \label{fcy}%
\end{equation}
where
\[
\mathcal{Y}_{L}\subset\mathbb{CP}^{N}\times\mathcal{H}_{L}.
\]
Let $\widetilde{\mathcal{H}}_{L}$ be the universal covering of $\mathcal{H}%
_{L}$ and let
\begin{equation}
\widetilde{\mathcal{Y}}_{L}\rightarrow\widetilde{\mathcal{H}}_{L} \label{fam1}%
\end{equation}
be the pullback family of $\left(  \ref{fcy}\right)  .$ It is easy to see that
the group $\mathbb{SL}_{N+1}(\mathbb{C})$ acts on $\mathcal{H}_{L}.$ This
implies that $\mathbb{SL}_{N+1}(\mathbb{C})$ acts also on $\widetilde
{\mathcal{H}}_{L}.$ We will need the following Lemma and its Corollary:

\begin{lemma}
\label{auto1}Let $G$ be a subgroup of $\mathbb{SL}_{N+1}(\mathbb{C})$ that
fixes a point $\tau_{0}\in\widetilde{\mathcal{H}}_{L},$ then $G$ is a finite
group of complex analytic automorphisms of the CY manifold M$_{\tau_{0}}$ and
it is a normal subgroup of $\mathbb{SL}_{N+1}(\mathbb{C})$ that acts trivially
on $H_{n}($M,$\mathbb{Z}).$
\end{lemma}

\textbf{Proof: }According to Theorem \ref{jt} $\mathcal{H}_{L}$ is a smooth
quasi-projective variety. So its universal cover $\widetilde{\mathcal{H}}_{L}$
is a simply connected complex manifold and we may suppose that the family
$\left(  \ref{fam1}\right)  $ is marked and polarized. The Definition
\ref{mark} of the marked family of CY manifolds implies that if $G$ fixes the
point $\tau_{0}\in\widetilde{\mathcal{H}}_{L},$ then $G$ is a subgroup in
$\mathbb{SL}_{N+1}(\mathbb{C})$ that stabilizes M$_{\tau_{0}}$ in
$\mathbb{CP}^{N}.$ This shows that $G$ must be a group of holomorphic
automorphisms of N$_{\tau_{0}}$ and it must act trivially on $H_{n}%
($M,$\mathbb{Z}).$ Theorem \ref{auto0} implies that $G$ acts trivially on
$\widetilde{\mathcal{H}}_{L}.$ This will imply that it is a normal subgroup in
$\mathbb{SL}_{N+1}(\mathbb{C}).$ To show this fact we need to prove that for
any element $g\in\mathbb{SL}_{N+1}(\mathbb{C})$ we have
\[
g^{-1}Gg=G.
\]
Let $\tau=g^{-1}(\tau_{0}).$ Then direct computations show that for any
element $h\in G$ we have
\[
g^{-1}hg(\tau)=g^{-1}\circ g(\tau)=\tau.
\]
The last equality shows that $G$ is a normal subgroup. Lemma \ref{auto1} is
proved. $\blacksquare.$

\begin{corollary}
\label{auto2}The group
\[
G_{1}:=\mathbb{SL}_{N+1}(\mathbb{C)}/G
\]
acts freely on $\widetilde{\mathcal{H}}_{L}.$
\end{corollary}

We will prove that the quotients $\widetilde{\mathcal{Y}}_{L}/G_{1}$ and
$\widetilde{\mathcal{H}}_{L}/G_{1}$ exist as complex manifolds and that
\[
\mathcal{U}_{L}:=\widetilde{\mathcal{Y}}_{L}/G_{1}\rightarrow\mathcal{T}%
_{L}(N):=\widetilde{\mathcal{H}}_{L}/G_{1}%
\]
will be the family $\left(  \ref{fam2}\right)  $ with the properties stated in
the Theorem.

Palais proved in \cite{Pa} the following Theorem:

\begin{theorem}
Suppose that an arbitrary Lie group $\mathfrak{G}$ acts on a complex manifold
Y, then the quotient Y/$\mathfrak{G}$ exists in the category of complex
spaces, provided that the map
\begin{equation}
\psi:\ \mathfrak{G}\times\text{Y}\rightarrow\text{Y}\times\text{Y }
\label{fam3}%
\end{equation}
defined as $\psi(g,m)=(g(m),m)$ is proper.
\end{theorem}

In order to prove that the morphism defined by the action of the Lie group
\[
\mathfrak{G}\mathfrak{=}G_{1}:=\mathbb{SL}_{k+1}(\mathbb{C})/G
\]
acts properly on $\widetilde{\mathcal{H}}_{L}$ and $\widetilde{\mathcal{Y}%
}_{L},$ we need to use the following criterion for the properness of the map
that can be found in \cite{EGA}, Chapter II, 7:

\begin{criterion}
\label{Groth}Let $\psi:\ $X$\rightarrow$Z be a morphism of algebraic varieties
over an algebraically closed field $k$, $\mathcal{O}$ a discrete valuation
ring with a residue class field $k,$%
\[
\phi^{\ast}:Spec\,\mathcal{O}\rightarrow\text{X}%
\]
is a rational morphism and%
\begin{equation}
h:\ Spec\,\mathcal{O}\rightarrow\text{Z} \label{cond0}%
\end{equation}
is a morphism, where
\begin{equation}
\psi\circ\phi^{\ast}=h. \label{cond}%
\end{equation}
If for such $\phi^{\ast}$ and h there exists a morphism
\[
\phi:\ Spec\,\mathcal{O\rightarrow}\text{X},
\]
coinciding with $\phi^{\ast}$ on the generic point of $Spec\,\mathcal{O},$ the
morphism f is proper.
\end{criterion}

To apply the Criterion \ref{Groth} in our situation we notice that the
properness of the map $\left(  \ref{fam3}\right)  $ for the manifolds
$\widetilde{\mathcal{Y}}_{L}$ and $\widetilde{\mathcal{H}_{L}}$ and the group
$G_{1}$ is an obvious formal consequence of the analogous fact for the
varieties $\mathcal{Y}_{L}$ and $\mathcal{H}_{L}.$ In this case we have to
deal with algebraic varieties and algebraic action of the group $G_{1}.$
Therefore $\psi$ is a morphism in the category of algebraic varieties.

The existence of the Grothendieck families $\mathcal{Y}_{L}\rightarrow
\mathcal{H}_{L}$ defines the following families
\begin{equation}
\pi_{1}:\mathcal{N}_{L}\rightarrow G_{1}\times\mathcal{H}_{L} \label{con}%
\end{equation}
and
\begin{equation}
\pi_{2}:\mathcal{R}_{L}\rightarrow\mathcal{H}_{L}\times\mathcal{H}_{L},
\label{cona}%
\end{equation}
where $\pi_{1}^{-1}(g,\tau):=$N$_{g(\tau)}\times$N$_{\tau}$ and $\pi_{2}%
^{-1}(\tau_{1},\tau_{2})=$N$_{\tau_{1}}\times$N$_{\tau_{2}}.$

Next we are going to apply Criterion \ref{Groth} to the map $\psi$ defined by
$\left(  \ref{fam3}\right)  $. We will need to use Theorem 2 of the paper
\cite{MM} of Matsusaka and Mumford. It asserts that if $V$ and $W$ are smooth
polarized varieties over $Spec\mathcal{O}$, neither of which is ruled, and
their specialization $V_{0}$ and $W_{0}$ over the closed point of
$Spec\mathcal{O}$ are also smooth polarized varieties, then the specialization
of any isomorphism $\psi:V\rightarrow W$ is an isomorphism $\psi_{0}%
:V_{0}\rightarrow W_{0}.$

It is easy to see that giving a morphism of algebraic varieties
\[
h:\ Spec\,\mathcal{O\rightarrow H}_{L}\times\mathcal{H}_{L}%
\]
and a rational morphism%
\[
\phi^{\ast}:\ Spec\,\mathcal{O}\rightarrow G_{1}\times\mathcal{H}_{L}%
\]
satisfying the criterion stated in the Grothendieck's Lemma, defines varieties
$V$ and $W$ over $Spec\,\mathcal{O}$ and the isomorphism $\psi^{\ast}$ of
them. On the other hand $\phi^{\ast}$ gives a rational morphism of
$Spec\,\mathcal{O}$ into $G_{1}$ and therefore a rational morphism $\phi
^{\ast}$ of $Spec\,\mathcal{O}$ into the projective closure of $\overline
{G}_{1}$ of the projective group $G_{1}.$ Since the dimension of the scheme
$Spec\,\mathcal{O}$ is one, then the rational morphism $\phi^{\ast}$ of
$Spec\,\mathcal{O}$ into a projective variety $\overline{G}_{1}$ can be
prolonged to a morphism $\phi$\ of $Spec\,\mathcal{O}$ to $\overline{G}_{1}.$
Thus we get a family $\overline{W}\rightarrow Spec\,\mathcal{O}$ which is
defined over the closed point of $Spec\,\mathcal{O}$ and contains the family
$W\rightarrow Spec\,\mathcal{O}.$ The fact that the rational morphism
$\phi^{\ast}$ can be prolonged to an algebraic morphism $\phi$ implies that
the morphism $\psi^{\ast}$ can be prolonged to a morphism $\psi$ between the
varieties $V$ and $\overline{W}.$ By the Theorem of Matsusaka and Mumford, the
specialization of $\psi$ is an isomorphism, i.e. $\psi$ takes closed point of
$Spec\,\mathcal{O}$ into a point of $G_{1}.$ This means that $\psi$ is an
algebraic morphism of $Spec\,\mathcal{O}$ into $G_{1}.$ Thus the criterion of
properness for the morphism $\psi$ holds. As we pointed out this fact implies
that the action of $G_{1}$ on $\widetilde{\mathcal{H}}_{L}$ is proper too. The
quotient $\widetilde{\mathcal{H}}_{L}/G_{1}$ exists and it is a smooth
algebraic variety.

The situation is analogous to the morphism
\[
\pi:\ \widetilde{\mathcal{Y}}_{L}\rightarrow\mathcal{U}_{L}=\widetilde
{\mathcal{Y}}_{L}/G_{1}.
\]
Thus we proved that the quotients $\widetilde{\mathcal{Y}}_{L}/G_{1}$ and
$\widetilde{\mathcal{H}}_{L}/G_{1}$ exist as complex spaces. It is obvious
that they define a family
\[
\mathcal{U}_{L}\rightarrow\mathcal{T}_{L}(\text{M})
\]
of marked polarized CY manifolds. According to Cor. \ref{auto2}, $G_{1}$ acts
on $\mathcal{Y}_{L}$ and $\mathcal{H}_{L}$ without fixed points. Therefore,
\begin{equation}
\widetilde{\mathcal{Y}}_{L}/G_{1}=\mathcal{U}_{L}\text{ and }\widetilde
{\mathcal{H}}_{L}/G_{1}=\mathcal{T}_{L}(\text{M}) \label{HilbI}%
\end{equation}
are manifolds. Thus we have constructed the family $\left(  \ref{fam2}\right)
$ as required in Theorem \ref{jt}$.$ We also proved that the base of the
family $\left(  \ref{fam2}\right)  $ $\mathcal{T}_{L}($M$)$ is smooth. Next we
will prove that the complex dimension of $\mathcal{T}_{L}($M$)$ is $h^{n-1,1}$
and that the tangent space at each point $\tau\in\mathcal{T}_{L}($N$)$ is
isomorphic to $H^{1}($M$_{\tau},\Theta_{\text{M}_{\tau}}).$

We shall prove that the family $\left(  \ref{fam2}\right)  $ is effectively
parametrized. According to Kodaira this means that for any point
$s\in\mathcal{T}($M) the mapping
\[
T_{s,\mathcal{T}_{L}(\text{M})}\rightarrow H^{1}(\text{M}_{s},\Theta
_{\text{M}_{s}})
\]
is an isomorphism. From the construction of $\mathcal{T}_{L}($M) as defined in
$\left(  \ref{HilbI}\right)  $ we can conclude that
\[
T_{s,\mathcal{T}_{L}(\text{M})}\approxeq H^{0}(\text{M},\mathcal{N}%
_{\mathbb{CP}^{N}/\text{M}})/sl(N+1)
\]
where $sl(N+1)$ is the Lie algebra of the group $\mathbb{SL(}N+1).$ In
\cite{JT96} we proved that
\[
H^{0}(\text{M},\mathcal{N}_{\mathbb{CP}^{N}/\text{M}})/sl(N+1)\approxeq
H^{1}(\text{M,}\Theta_{\text{M}}).
\]
Theorem \ref{jt} is proved. $\blacksquare$

\begin{corollary}
\label{teich1}Let $\mathcal{Y\rightarrow}$X be any family of marked CY
manifolds, then there exists a unique holomorphic map
\[
\phi:\ \text{X}\rightarrow\mathcal{T}_{L}(M)
\]
up to a biholomorphic map $\psi$ of M which induces the identity map on
$H_{n}($M$,\mathbb{Z}).$
\end{corollary}

From now on we will denote by $\mathcal{T}$(M) the irreducible component of
the Teichm\"{u}ller space that contains our fixed CY manifold M.

\subsection{The Mapping Class Group}

\begin{definition}
\label{mcg}We will define the mapping class group $\Gamma($M) of any compact
C$^{\infty}$ manifold M as follows: $\Gamma=Diff_{+}($M)$/Diff_{0}($M)$,$
where $Diff_{+}($M$)$ is the group of diffeomorphisms of M preserving the
orientation of M and $Diff_{0}($M$)$ is the group of diffeomorphisms isotopic
to identity.
\end{definition}

\begin{definition}
Let $L\in H^{2}($M$,\mathbb{Z})$ be the imaginary part of a K\"{a}hler metric.
Let
\[
\Gamma_{L}:=\{\phi\in\Gamma\text{(M)}|\phi(L)=L\}.
\]

\end{definition}

It is a well know fact that the moduli space of polarized algebraic manifolds
$\mathfrak{M}_{L}$(M)$=\mathcal{T}_{L}($M$)/\Gamma_{L}.$

\begin{theorem}
\label{Vie}There exists a subgroup of finite index $\Gamma$ of $\ \Gamma_{L}$
such that $\Gamma$ acts freely on $\mathcal{T}$(M) and $\Gamma\backslash
\mathcal{T}_{L}$(M)$=\mathcal{M}_{L}$(M) is a non-singular quasi-projective variety.
\end{theorem}

\textbf{Proof: }In case of odd dimensional CY manifolds there is a
homomorphism induced by the action of the diffeomorphism group on the middle
homology with coefficients in $\mathbb{Z}:$%
\[
\phi:\Gamma_{L}\rightarrow\mathbb{S}p(2b_{n},\mathbb{Z)}.
\]
In the case of even dimensional CY, there is a homomorphism
\[
\phi:\Gamma_{L}\rightarrow\mathbb{SO}(2p,q;\mathbb{Z)}%
\]
where $\mathbb{SO}$($2p,q;\mathbb{Z})$ is the group of the automorphisms of
the lattice $H_{n}$(M,$\mathbb{Z})$/Tor. An important theorem due to Sullivan
proved in \cite{Sul} states:

\begin{theorem}
Suppose that the real dimension of a C$^{\infty}$ manifold M is bigger or
equal to 5, then the image $\phi(\Gamma_{L})$ of the mapping class group is an
arithmetic group.
\end{theorem}

This result of Sullivan implies that $\Gamma_{L}$ is an arithmetic group. So
the image of $\Gamma_{L}$ has a finite index in the groups $\mathbb{S}%
p$($2b_{n},\mathbb{Z)},$ $\mathbb{SO}(2p,q;\mathbb{Z)}$\ and so $\ker(\phi)$
is a finite group. A theorem of Borel implies that we can always find a
subgroup of finite index $\Gamma$ in $\Gamma_{L}$ such that $\Gamma$ acts
freely on $\mathbb{S}p$($2b_{n},\mathbb{R)}$/$\mathbb{U}$($b_{n}\mathbb{)}$ or
on $\mathbb{SO}_{0}$($2p,q;\mathbb{R})$/$\mathbb{SO}$($2p$)$\times\mathbb{SO}%
$($q$). We will prove that $\Gamma$ acts without fixed point on $\mathcal{T}%
_{L}$(M).

Let $\mathcal{K}$ be the Kuranishi space of the deformations of M. Suppose
that there exists an element $g\in\Gamma$, such that $g(\tau)=\tau$ for some
\[
\tau\in\mathcal{K\subset T}_{L}(\text{M}).
\]
From the local Torelli theorem we deduce that we may assume that the Kuranishi
space $\mathcal{K}$ is embedded in $\mathcal{G}$, the classifying space of the
Hodge structures of weight n on $H^{n}($M,$\mathbb{Z)\otimes C}$. Griffiths
proved in \cite{Gr} that $\mathcal{G\thickapprox}G/K$ \ where G in the odd
dimensional case is $\mathbb{S}p$($2b_{n},\mathbb{R)}$ and in the even
dimensional is $\mathbb{SO}_{0}$($2p,q;\mathbb{R})$ and $K$ is a compact
subgroup of $G.$

Let $K_{0}$ be the maximal compact subgroup of $G.$ So we have a natural
$C^{\infty}$ fibration
\[
K_{0}/K\subset G/K\rightarrow G/K_{0}.
\]
Griffith's transversality theorem implies that $\mathcal{K}$ is transversal to
the fibres $K_{0}/K$ of the fibration $G/K\rightarrow G/K_{0}$.

The first part of our theorem follows from the fact that $\mathcal{K}$ is
transversal to the fibres $K_{0}/K$ of the fibration $G/K\rightarrow G/K_{0}$
and the following observation; if $g\in\Gamma$ fixes a point $\tau\in G/K_{0}
$, then $g\in K_{0}\cap\Gamma$.\footnote{We suppose that $K$ \ or $K_{0}$ acts
on the right on $G$ and $\Gamma$ acts on the left on $G.$} On the other hand
side it is easy to see that the local Torelli theorem implies that the action
of $\Gamma$ on $\mathcal{K}$ is induced from the action $\Gamma$ on $G/K$ by
left multiplications. So we can conclude that the action of $\Gamma$ preserves
the fibration
\[
K_{0}/K\subset G/K\rightarrow G/K_{0}.
\]
The first part of our theorem follows directly from here and the fact that
$\Gamma$ acts without fix point on $G/K_{0}$.

The second part of the theorem, namely that the space $\Gamma$%
$\backslash$%
$\mathcal{T}_{L}$(M) is a quasi projective follows directly from the fact that
$\Gamma$%
$\backslash$%
$\mathcal{T}_{L}$(M)$\rightarrow\Gamma_{L}$%
$\backslash$%
$\mathcal{T}_{L}$(M) is a finite map and that $\Gamma_{L}$%
$\backslash$%
$\mathcal{T}_{L}$(M) is a quasi projective variety according to \cite{V}.
Theorem \ref{Vie} is proved. $\blacksquare$

\subsection{Construction of the Moduli Space of Polarized CY Manifolds}

According to Viehweg the coarse moduli space of polarized CY manifolds
\[
\mathfrak{M}_{L}(\text{M})=\mathcal{H}_{L}/\mathbb{SL}_{N+1}(\mathbb{C})
\]
is a quasi-projective variety. See \cite{V}. On the other hand it is a
standard that the coarse moduli space is just the following quotient:
\[
\mathfrak{M}_{L}(\text{M})=\mathcal{T}_{L}(\text{M})/\Gamma_{L}.
\]

\begin{theorem}
\label{Mod}There exists a finite cover $\mathcal{M}_{L}($M$)$ of
\[
\mathfrak{M}_{L}(\text{M})=\mathcal{T}_{L}(\text{M})/\Gamma_{L}%
\]
with the following properties: \textbf{A.} $\mathcal{M}_{L}($M$)$ is a smooth
algebraic variety, \textbf{B. }Over $\mathcal{M}_{L}($M$)$ there exists a
family
\begin{equation}
\mathcal{N}_{L}\rightarrow\mathcal{M}_{L}(\text{M}) \label{unf}%
\end{equation}
of polarized CY manifolds with the following property; Let
\begin{equation}
\pi_{\mathcal{C}}:\mathcal{Y\rightarrow C} \label{qu}%
\end{equation}
be any complex analytic family of polarized CY manifolds with a class of
polarization $L$ such that at least one of its fibres of the family $\left(
\ref{qu}\right)  $ is isomorphic as a polarized variety to a fibre of the
family $\left(  \ref{unf}\right)  $. Then there exists a unique complex
analytic map%
\[
\phi:\mathcal{C\rightarrow M}_{L}(\text{M})\text{.}%
\]
which induces a holomorphic map between the families $\left(  \ref{qu}\right)
$and $\left(  \ref{unf}\right)  .$ The map between the families is defined
uniquely up to a biholomorphic map $\phi$ of M which induces the identity map
on $H_{n}($M$,\mathbb{Z}).$
\end{theorem}

\textbf{Proof: }Let
\begin{equation}
\psi:\Gamma_{L}\rightarrow H^{n}(\text{M,}\mathbb{Z}) \label{Triv}%
\end{equation}
be the natural representation of the group $\Gamma_{L}.$ The results of
Sullivan imply that the image $\psi(\Gamma_{L})$ is an arithmetic group.
According to a Theorem of \ref{Vie} there exists a subgroup $\Gamma_{L}^{{"}}$
of finite index in $\psi(\Gamma_{L})$ such that $\psi(\Gamma_{L})$ acts freely
on the Teichm\"{u}ller space $\mathcal{T}_{L}$(M). From here we deduce that
the space $\Gamma_{L}^{{"}}\backslash\mathcal{T}_{L}$(M) is smooth. Let us
define $\Gamma_{L}^{{\prime}}:=\psi^{-1}(\Gamma_{L}^{{"}}).$ Clearly
$\Gamma_{L}^{{\prime}}$ is a subgroup of finite index in $\Gamma_{L}.$ The
local Torelli theorem implies that $\mathcal{T}_{L}($M)$/\Gamma_{L}^{{\prime}%
}$ will be a non-singular variety. From the definition of $\psi$ given by
$\left(  \ref{Triv}\right)  $ it follows that $\ker\psi=G$ acts trivially on
$H^{n}($M,$\mathbb{Z}).$ The local Torelli Theorem implies that $G$ acts
trivially on the Teichm\"{u}ller space $\mathcal{T}_{L}($M)$.$ The existence
and the properties of the family%
\[
\mathcal{N}_{L}\rightarrow\mathcal{M}_{L}(\text{M})
\]
follows from the existence and the properties of the family
\[
\mathcal{Z}_{L}\rightarrow\mathcal{T}_{L}(\text{M})
\]
proved in Theorem \ref{teich1}. $\blacksquare$

\subsection{Compactifications of the Moduli Spaces $\mathcal{M}_{L}($M$)$ and
$\mathfrak{M}_{L}($M)}

According to \cite{V} the coarse moduli space $\mathfrak{M}_{L}($M) is a
quasi-projective variety. Since $\mathcal{M}_{L}($M$)$ is a finite cover of
$\mathfrak{M}_{L}($M), then $\mathcal{M}_{L}($M) will be quasi-projective too.
Hironaka's Theorem about resolution of singularities states that we can always
compactify both $\mathcal{M}_{L}($M$)$ and $\mathfrak{M}_{L}($M) such that
\begin{equation}
\overline{\mathcal{M}_{L}(\text{M})}\ominus\mathcal{M}_{L}(\text{M}%
)=\mathcal{D}_{\infty}\text{ }and\text{ }\overline{\mathfrak{M}_{L}(\text{M}%
)}\ominus\mathfrak{M}_{L}(\text{M})=\mathfrak{D}_{\infty} \label{vie}%
\end{equation}
the divisors $\mathcal{D}_{\infty}$ and $\mathfrak{D}_{\infty}$ are divisors
with normal crossings.

\section{Moduli of Maps}

\subsection{Basic Definitions}

Let C and X be two projective varieties. Let $f:$ C$\rightarrow$X be a
morphism between them. Let $\Gamma_{f}\subset$C$\times$X be the graph of the
map $f:$C$\rightarrow$X$.$ According to the results of Grothendieck, the
Hilbert scheme of $\Gamma_{f}\subset$C$\times$X is a projective scheme. See
\cite{EGA}.

\begin{definition}
\label{Gro}We will denote the Hilbert scheme of $\Gamma_{f}\subset$C$\times$X
by $\mathfrak{M}_{f}($C$,$X$)$ and we will call it the moduli space of the map
$f.$
\end{definition}

The above mentioned results of Grothendieck implies that $\mathfrak{M}_{f}%
($C$,$X$)$ is a projective scheme.

\begin{definition}
\label{d1}\textbf{a.} Let f: C$\rightarrow$X be a morphism of projective
varieties with
\[
dim_{\mathbb{C}}\text{C}\leq dim_{\mathbb{C}}\text{X}.
\]
Suppose that the morphism
\[
\text{f}:\text{C}\rightarrow\text{f}(\text{C})
\]
is finite. We say that f admits a non-trivial one parameter deformations if
there is a non-singular projective curve T and a family of algebraic maps
\[
\text{F}:\text{T}\times\text{C}\rightarrow\text{X}%
\]
such that for some t$_{0}\in$T we have F$_{t_{0}}=$f and the morphism
\[
\text{F}:\text{T}\times\text{C}\rightarrow\text{F}(\text{T}\times\text{C})
\]
is finite too. \textbf{b. }\textit{We will say that the deformation of f is
trivial if }
\[
\text{\textit{F}}_{t}=\text{\textit{f}}%
\]
\textit{\ for all t}$\in$T. \textbf{c.}\textit{\ Two families of maps }
\[
\text{\textit{F}}_{1}:\text{T}\times\text{\textit{C}}_{1}\rightarrow
\text{\textit{X}}%
\]
\textit{and }
\[
\text{\textit{F}}_{2}:\text{T}\times\text{\textit{C}}_{2}\rightarrow
\text{\textit{X}}%
\]
\textit{\ are said to be isomorphic if there is a common finite cover C of
C}$_{1}$\textit{\ and C}$_{2}$\textit{\ such that the lifts of F}$_{1}%
$\textit{\ to T}$\times$\textit{C and F}$_{2}$ \textit{to T}$\times$\textit{C
are isomorphic, meaning there exists a biholomorphic map id}$\times$\textit{g
from T}$\times$\textit{C to itself such that}
\[
\text{\textit{F}}_{1}{}=\text{F}_{2}\circ(\text{id}\times\text{g}).\mathit{\ }%
\]

\end{definition}

From now on we will consider only one parameter deformations of maps.

\begin{definition}
\label{Sh}Let C be a fixed Riemann surface. Let X be a quasi-projective space
and $\overline{\text{X}}$ some projective compactification of X such that
$\mathcal{D}_{\infty}=\overline{\text{X}}$ $\circleddash$ X is a divisor with
normal crossings. Assume that on X there exists a K\"{a}hler metric with
non-positive holomorphic sectional curvature bounded away from zero and
logarithmic growth near $\mathcal{D}_{\infty}$. We will consider from now on
the set $\mathcal{M}_{\text{C}}$ of all holomorphic maps
\[
\text{f}:\text{C}\rightarrow\overline{\text{X}}%
\]
such that: \textbf{a. } the map f:C$\rightarrow$f(C) is finite map and
\textbf{b.} f(C) \textit{is not contained in\ }$\mathcal{D}_{\infty}$ and
\textit{the pullback }
\[
\text{f}^{\ast}(\mathcal{D}_{\infty})=\text{D}_{\text{S}},
\]
\textit{where }D$_{\text{S}}$ \textit{is a fixed divisor on C. }
\end{definition}

\subsection{Smoothness of the Moduli Space of Maps}

\begin{theorem}
\label{kefeng}Suppose that X is a projective manifold such that on a Zariski
open set
\[
\text{X}^{{\prime}}:=\overline{\text{X}}\ominus\mathcal{D}_{\infty}%
\]
there exists a K\"{a}hler metric with a non-positive holomorphic bi-sectional
curvature. Let C be a fixed algebraic curve. Let
\begin{equation}
\mathcal{M}_{\text{C}}:\left\{  \text{f}:\text{C}\rightarrow\text{X%
$\vert$%
f is a finite hol map from C to f(C) and f(C)}\varsubsetneq\mathcal{D}%
_{\infty}\right\}  . \label{kefeng0}%
\end{equation}
Then the moduli space $\mathcal{M}_{\text{C}}$ is a smooth quasi-projective variety.
\end{theorem}

\textbf{Proof: }The proof is based on several lemmas. The following lemma
describes the tangent space to the moduli space $\mathcal{M}_{\text{C}}.$

\begin{lemma}
\label{0l}Let
\[
\text{f}:\text{C}\rightarrow\text{X}%
\]
be a finite map from an algebraic curve C to f(C) such that
\[
\text{f}(\text{C})\varsubsetneq\mathcal{D}_{\infty}.
\]
Let
\[
\sigma\in H^{0}(\text{C},f^{\ast}(T_{\text{X}}))
\]
be a non-zero section with a bounded norm. Then $\sigma$ is parallel with
respect to the connection induced by the restriction of the K\"{a}hler metric
on X$\ominus\mathcal{D}_{\infty}$.
\end{lemma}

\textbf{Proof: }Denote by $\text{C}^{0}=\text{C}\ominus$f$^{\ast}%
(\mathfrak{D}_{\infty})$. Let
\[
\sigma\in H^{0}(\text{C}\ominus\text{f}^{\ast}(\text{f}_{\ast}\text{(C)}%
\cap\mathcal{D}_{\infty}),f^{\ast}(T_{\text{X}}))
\]
be a non-zero section. Suppose that
\begin{equation}
\nabla\sigma\neq0, \label{k0}%
\end{equation}
where $\nabla$ is the covariant derivative on the vector bundle f$^{\ast}%
$(T$_{\text{X}})|_{\text{C}^{0}}$ induced by the pullback of the K\"{a}hler
metric on X. Direct computations show that the following formula holds:%

\begin{equation}
\frac{d^{2}}{d\tau\overline{d\tau}}\left\Vert \sigma\right\Vert ^{2}%
=\left\Vert \nabla\sigma\right\Vert ^{2}-\left\langle R({\frac{d}{d\tau
},\overline{\frac{d}{d\tau}}})\sigma,\sigma\right\rangle \geq0, \label{k1}%
\end{equation}
where $\left\langle R({\frac{d}{d\tau},\overline{\frac{d}{d\tau}}}%
)\sigma,\sigma\right\rangle $ is the corresponding holomorphic bi-sectional
curvature on the vector bundle f$^{\ast}$(T$_{\text{X}})|_{\text{C}^{0}}$ with
a metric induced by the pullback of the K\"{a}hler metric on X. Formula
$\left(  \ref{k1}\right)  $ and the assumption that we can prolong the finite
map f from $\text{C}^{0}=\text{C}\ominus$f$^{\ast}(\mathfrak{D}_{\infty})$ to
a projective map show that the function $\left\Vert \sigma\right\Vert ^{2}$ is
a bounded plurisubharmonic function on $\text{C}^{0}$. Therefore by maximal
principle it is the constant function because formula $\left(  \ref{k1}%
\right)  $ implies
\begin{equation}
\frac{d^{2}}{d\tau\overline{d\tau}}\left\Vert \sigma\right\Vert ^{2}=0.
\label{k1a}%
\end{equation}
So we get a contradiction. This implies that $\sigma=0$ if $\sigma$ is not
parallel. So if $\sigma$ is a holomorphic non-zero section with finite $L^{2}$
norm, then it must be parallel. Lemma \ref{0l} is proved. $\blacksquare$

\begin{corollary}
\label{00l}The tangent space at the point f$\in\mathcal{M}_{\text{C}}$ is
isomorphic to all parallel sections $\sigma\in H^{0}($C$,f^{\ast}(T_{\text{X}%
})).$
\end{corollary}

\begin{lemma}
\label{-1l}Let
\[
\sigma\in H^{0}(\text{C}\ominus\text{f}^{\ast}(\text{f(C)}\cap\mathcal{D}%
_{\infty}),f^{\ast}(T_{\text{X}})).
\]
be a parallel section. Then there exists a family of maps
\[
\text{F}_{t}:\ \text{C}\rightarrow\text{X}%
\]
where $t\in\mathcal{D}$ and $\mathcal{D}$ is the unit disk such that
\[
\text{F}_{0}=\text{f}%
\]
and
\[
\frac{d}{dt}\text{F}_{t}|_{t=0}=\sigma.
\]

\end{lemma}

\textbf{Proof: }Since
\[
\sigma\in H^{0}(\text{C}\ominus\text{f}^{\ast}(\text{f(C)}\cap\mathcal{D}%
_{\infty}),\text{f}^{\ast}(T_{\text{X}}))
\]
is a parallel section, there exists a vector $\overrightarrow{\tau}\in
T_{x,\text{X}}$ at each point $x\in$f$($C$)\ominus\mathcal{D}_{\infty}$ such
that
\[
\text{f}^{\ast}(\overrightarrow{\tau})=\sigma
\]
and $x\in$f(C). Since the metric $g$ on X is K\"{a}hler we can define as in
the Paragraph \ref{cexp}, Cor. \ref{lines} the exponential complex map:
\[
\exp_{x}:T_{x,\text{X}}\rightarrow\text{X}%
\]
at each point $x\in$f(C) such that
\begin{equation}
\exp_{x}(t\overrightarrow{\tau})=D_{x} \label{exp0}%
\end{equation}
is a totally geodesic disk when $|t|<1$ and the tangent vector of the disk at
the point $x$ is $\overrightarrow{\tau}.$ Here we are using the fact that the
holomorphic bi-sectional curvature is non-positive. This implies that the
complex exponential map is injective for $|t|<1.$ Using the fact that $\sigma$
is a parallel section of the vector bundle f$^{\ast}T_{\text{X}}$ on
C$\ominus$f$^{\ast}($f(C)$\cap\mathcal{D}_{\infty})$ we deduce immediately
that by defining
\begin{equation}
\text{f}_{t}(\text{s})=\exp_{\text{f(s)}}(t\overrightarrow{\tau}) \label{exp1}%
\end{equation}
and taking the exponential map $\left(  \ref{exp0}\right)  $ at each point
$x\in f$(C), we have constructed a one parameter deformation
\[
\text{f}_{t}:\text{C}\times\mathcal{D}\rightarrow\text{X}%
\]
of f: C$\rightarrow$X as required in lemma \ref{-1l}. $\blacksquare.$

Lemma \ref{-1l} and Corollary \ref{00l} imply that the moduli space of all the
maps
\[
\text{f}:\text{C}\rightarrow\text{X}%
\]
such that f(C)$\subsetneq\mathcal{D}_{\infty}$ is smooth. Indeed we know that
the tangent space $T_{\text{f,}\mathcal{M}_{\text{C}}}$ of the moduli space
$\mathcal{M}_{\text{C}}$ at f$\in\mathcal{M}_{\text{C}}$ is isomorphic to
$H^{0}($C,f*($T_{\text{X}})).$ Cor. \ref{00l} tells us that $T_{\text{f,}%
\mathcal{M}_{\text{C}}}$ consists of all parallel sections. From Lemma
\ref{-1l} it follows that each parallel section can be integrated to a family
of maps parametrized by the unit disk. Grothendieck's theory of Hilbert
schemes tells us that the Hilbert scheme of all maps
\[
\text{f}:\text{C}\rightarrow\text{X}%
\]
such that
\[
\text{f}(\text{C})\nsubseteq\mathcal{D}_{\infty}%
\]
is a quasi-projective space. Theorem \ref{kefeng} is proved. $\blacksquare$

\subsection{Some Applications}

Let $($C$;x_{1},..,x_{n})$ be a Riemann surface with n distinct fixed points
on it. Let
\[
\mathcal{X}_{\text{C}}\rightarrow\text{C}%
\]
be a family of polarized CY manifolds with degenerate fibres over the points
$x_{1},..,x_{n}$. It is easy to see that the family $\mathcal{X}\rightarrow$C
defines a map
\begin{equation}
f:\text{C}\circleddash(x_{1}\cup...\cup x_{n})\rightarrow\mathcal{T}%
_{L}(\text{M})/\Gamma_{L}. \label{F}%
\end{equation}

\begin{lemma}
\label{TOD0}The map $f$ can be prolonged to a map
\[
\overline{f}:\text{C}\rightarrow\overline{\Gamma_{L}\backslash\mathcal{T}%
_{L}(\text{M})}=\overline{\mathfrak{M}_{L}(\text{M})}.
\]
where the compactification $\overline{\mathfrak{M}_{L}(\text{M})}$ of
$\mathfrak{M}_{L}($M$)$ is defined in $\left(  \ref{vie}\right)  $.
\end{lemma}

\textbf{Proof: }The proof of Lemma \ref{TOD0} is exactly Proposition
\textbf{9.10 }proved in \cite{Gr}. Lemma \ref{TOD1a} is proved. $\blacksquare
.$

Let $\Gamma_{L}$ and $\Gamma$ be the arithmetic groups defined in Theorem
\ref{Vie}. The map
\[
\pi:\mathcal{T}_{L}(\text{M})/\Gamma=\mathfrak{M}_{L}(\text{M})\rightarrow
\mathcal{T}_{L}(\text{M})/\Gamma_{L}=\mathcal{M}_{L}(\text{M)}%
\]
is a finite morphism since $\Gamma_{L}/\Gamma$ is a finite set and its
cardinality is $N$. We will prove the following lemma:

\begin{lemma}
\label{TOD1}There exists a finite affine cover
\[
\psi:\text{C}_{1}\rightarrow\text{C}\ominus(x_{1}\cup...\cup x_{n})
\]
of degree $N=\#\Gamma_{L}/\Gamma$ such that the map $f$ defined by $\left(
\ref{F}\right)  $ can be lifted to a map
\begin{equation}
f_{1}:\text{C}_{1}\rightarrow\mathcal{T}_{L}(\text{M})/\Gamma. \label{F1}%
\end{equation}

\end{lemma}

\textbf{Proof: }The construction of the affine Riemann surface C$_{1}$ is done
in a standard way, namely
\[
\text{C}_{1}=\left(  \text{C}\ominus(x_{1}\cup...\cup x_{n})\right)
\times_{f\left(  \text{C}\ominus(x_{1}\cup...\cup x_{n})\right)  }\pi
^{-1}(f\left(  \text{C}\ominus(x_{1}\cup...\cup x_{n})\right)  .
\]
We know from the local Torelli theorem and the fact $\Gamma_{L}\backslash
\mathcal{T}_{L}($M$)$ is the moduli space of polarized CY manifolds, that the
map $f$ defined by $\left(  \ref{F}\right)  $ exists. We know that, since
$\Gamma\subset\Gamma_{L}$ is a subgroup of finite index, the map
\[
\pi:\Gamma_{L}\backslash\mathcal{T}_{L}(\text{M})\rightarrow\Gamma
\backslash\mathcal{T}_{L}(\text{M})
\]
is a finite map. Therefore the map $f$ can be lifted to a map
\[
f^{\prime}:\text{C}_{1}\rightarrow\mathcal{T}_{L}(\text{M})/\Gamma
_{L}(\text{M}).
\]
$\,$follows directly from the construction of the open Riemann surface
C$^{\prime}.$ Lemma \ref{TOD1} is proved. $\blacksquare$

\begin{lemma}
\label{TOD1a}The map $f_{1}$ constructed in Lemma \ref{TOD1} can be prolonged
to a map
\[
f:\overline{\text{C}_{1}}\rightarrow\overline{\Gamma\backslash\mathcal{T}%
_{L}(\text{M})}=\overline{\mathcal{M}_{L}(\text{M})}%
\]
where $\overline{\text{C}_{1}}$ is the \textit{closure of }C$_{1}$ in the
compactification $\overline{\mathcal{M}_{L}(\text{M})}$ of $\mathcal{M}_{L}%
($M$)$ as defined in $\left(  \ref{vie}\right)  $.
\end{lemma}

\textbf{Proof:} The proof of Lemma \ref{TOD1a} is the same as Lemma \ref{TOD0}
and is due to Griffiths. See \cite{Gr}. $\blacksquare$

\begin{definition}
\label{Rigid}We will say that a family of polarized CY manifolds
\[
\pi:\mathcal{X}_{\text{C}}\rightarrow\text{C}%
\]
over a Riemann surface C is rigid if the moduli space of the map
\[
\text{f}:\text{C}\rightarrow\mathfrak{M}_{L}\text{(M) }%
\]
induced by the natural holomorphic map of the base of the family into the
moduli space of polarized CY manifolds is a discrete set.
\end{definition}

The proof of the following remark is obvious and we will omit it.

Let $\mathcal{X\rightarrow}$C be a family of polarized CY manifolds over a
fixed Riemann surface C with a fixed points of degenerations. Let $\pi
_{1}:\overline{\text{C}_{1}}\rightarrow$C be the finite cover constructed in
Lemmas \ref{TOD1} and \ref{TOD1a}. Let $\mathcal{X}_{1}\rightarrow$
$\overline{\text{C}_{1}}$ is the pullback of the family $\mathcal{X\rightarrow
}$C. Then the family $\mathcal{X\rightarrow}$C is rigid if and only if the
pullback family $\mathcal{X}_{1}\rightarrow$ $\overline{\text{C}_{1}}$ is rigid.

\section{Yukawa Coupling and Rigidity}

\subsection{Yukawa Coupling}

\begin{definition}
\label{yuk}Suppose that $\omega_{\tau}$ is a family of holomorphic forms on
the local universal family of CY manifolds over the Kuranishi space of
polarized CY manifold. The expression
\[
\left\langle \underset{n}{\underbrace{\left(  \nabla_{\frac{\partial}%
{\partial\tau^{i_{1}}}}\circ\left(  \nabla_{\frac{\partial}{\partial
\tau^{i_{2}}}}\circ...\circ\left(  \nabla_{\frac{\partial}{\partial\tau
^{i_{n}}}}\left(  \omega_{\tau}\right)  \right)  \right)  \right)  }}%
,\omega_{\tau}\right\rangle
\]
is a holomorphic section of \ $\left(  \left(  \Omega_{\mathcal{X}%
/\mathcal{M}\text{(M)}}^{n}\right)  ^{\ast}\right)  ^{\otimes2}\otimes
S^{\otimes n}(\mathcal{T}_{\mathcal{M}\text{(M}}).$ This section is called the
Yukawa coupling.
\end{definition}

\begin{condition}
\label{Irr0}Suppose that for every $\phi\in H^{1}($M,$\Theta_{\text{M}})$
different from zero the following condition holds:
\begin{equation}
\wedge^{n}\phi\neq0\in{H}^{n}(\text{M},\wedge^{n}\Theta_{\text{M}})
\label{Irr}%
\end{equation}
Then we will prove that any family
\[
\pi:\mathcal{X}_{\text{C}}\rightarrow\text{C}%
\]
over a Riemann surface C is rigid if one of its fibre is M.
\end{condition}

The condition $\ref{Irr0}$ is difficult to check. Now we will formulate an
equivalent condition which is much easier to verify. Let
\begin{equation}
\pi:\mathcal{X}\rightarrow\mathcal{D} \label{fam}%
\end{equation}
be a non-trivial family of non singular n complex dimensional CY manifolds
over the unit disk $\mathcal{D}.$ We will denote by $t$ the local coordinate
in $\mathcal{D}.$ Let us denote by $\nabla_{\frac{\partial}{\partial t}}$ the
covariant differentiation induced by the Gauss-Manin connection. According to
Kodaira-Spencer-Kuranishi theory the tangent vectors at $0\in\mathcal{D}$ can
be identified with non zero elements $\phi\in H^{1}($M$_{0},\Theta
_{\text{M}_{0}}).$ Let $\Omega_{\mathcal{X}/\mathcal{D}}^{n}$ be the relative
dualizing sheaf. Let
\[
\omega_{t}\in H^{0}(\mathcal{X},\Omega_{\mathcal{X}/\mathcal{D}}^{n})
\]
be a family of holomorphic $n$ forms $\omega_{t}.$ We may suppose that
$\omega_{0}\neq0.$

\begin{theorem}
\label{Betty1} Let $\phi\in H^{1}($M,$\Theta_{\text{M}}),$ then the following
condition
\begin{equation}
\left\langle \underset{n}{\underbrace{\left(  \nabla_{\frac{\partial}{\partial
t}}\circ\left(  \nabla_{\frac{\partial}{\partial t}}\circ...\circ\left(
\nabla_{\frac{\partial}{\partial t}}\left(  \omega_{t}\right)  \right)
\right)  \right)  }}|_{t=0},\omega_{0}\right\rangle \neq0, \label{Be1}%
\end{equation}
is equivalent to $\left(  \ref{Irr}\right)  $, i.e. $\wedge^{n}\phi\neq0$ in
$H^{n}($M,$\wedge^{n}\Theta_{\text{M}}).$
\end{theorem}

\textbf{Proof:} It is a standard fact that for any family of forms $\omega
_{s}$ the covariant differentiation given by the Gauss-Manin connection is
given by the formula:
\begin{equation}
\nabla_{\frac{\partial}{\partial t}}\left(  \omega_{t}\right)  |_{t=0}%
=\omega_{0}\lrcorner\phi\label{Be2}%
\end{equation}
where $\lrcorner$ means contraction of tensors. From $\left(  \ref{Be2}%
\right)  $ we derive directly:
\begin{equation}
\underset{n}{\underbrace{\left(  \nabla_{\frac{\partial}{\partial t}}%
\circ\left(  \nabla_{\frac{\partial}{\partial t}}\circ...\circ\left(
\nabla_{\frac{\partial}{\partial t}}\left(  \omega_{t}\right)  \right)
\right)  \right)  }}|_{t=0}=\left(  \wedge^{n}\phi\right)  \lrcorner\omega
_{0}. \label{Be3}%
\end{equation}
$\left(  \ref{Be3}\right)  $ means that the form $\left(  \wedge^{n}%
\phi\right)  \lrcorner\omega_{0}$ is of type $(0,n).$ This implies directly
that $\left(  \ref{Be1}\right)  $ and $\left(  \ref{Irr}\right)  $ are
equivalent. Theorem \ref{Betty1} is proved. $\blacksquare$

\begin{theorem}
\label{TOD3}Let
\begin{equation}
\pi:\mathcal{X}_{\text{C}}\mathcal{\rightarrow}\text{C} \label{Tod1}%
\end{equation}
be a non-isotrivial family of CY manifolds over a compact Riemann surface C.
Suppose that there exists a point $\tau_{0}\in$C such that the fibre $\pi
^{-1}(\tau_{0})=$M$_{\tau_{0}}$ satisfies Condition $\ref{Irr0}$. Suppose that
$\mathcal{X}_{\text{C}}$ is a projective manifold. Then the family $\left(
\ref{Tod1}\right)  $ is rigid.
\end{theorem}

\textbf{Proof: }Theorem \ref{TOD3} follows from the fact that the moduli space
of the deformations of the map f defined in Definition \ref{Rigid} is smooth
by Theorem \ref{kefeng} and Lemma \ref{TOD31}. The rigidity therefore implies
that the tangent space of the moduli space at the point
\[
\text{f:\ C}\rightarrow\text{ }\overline{\mathcal{M}_{L}(\text{M})}%
\]
is zero dimensional.

\begin{lemma}
\label{TOD31}Suppose that the family $\left(  \ref{Tod1}\right)  $ satisfies
the conditions of Theorem \ref{TOD3}, then the set of the parallel sections of
the bundle $f^{\ast}(T_{\mathcal{M}_{L}\text{(M)}})$ consists only of the zero section.
\end{lemma}

\textbf{Proof: }The proof of Lemma \ref{TOD31} is done by contradiction.
According to \cite{Pe} the non rigidity of the family $\left(  \ref{Tod1}%
\right)  $ means the existence of
\[
\phi\in H^{0}(\text{C}\ominus\text{S},f^{\ast}(T_{\mathcal{M}\text{(M)}%
})),\text{ }\phi\neq0
\]
such that $\phi$ is parallel section with respect to the pullback of the Hodge
metric on $\mathcal{M}($M$)$. Since the tangent space at a point $\tau
\in\mathcal{M}($M$)$ can be identified with $H^{1}($M$_{\tau},T_{\text{M}%
}^{1,0}),$ the restriction of the section $\phi$ at each point $\tau\in
$C$\ominus$S can be viewed as a Kodaira-Spencer class
\[
\phi_{\tau}\in H^{1}(\text{M}_{\tau},T_{\text{M}}^{1,0}),
\]
Locally on M$_{\tau}$ $\phi_{\tau}$ is given by
\[
\phi_{\tau}:=\sum_{i,j}\left(  \phi_{\tau}\right)  _{\overline{j}}%
^{i}\overline{dz^{j}}\otimes\frac{d}{dz^{i}}.
\]
On the other hand $\wedge^{n}\phi_{\tau}$ is represented on M$_{\tau}$ by
\begin{equation}
\wedge^{n}\phi_{\tau}=\det\left(  \left(  \phi_{\tau}\right)  _{\overline{j}%
}^{i}\right)  \overline{dz^{1}}\wedge...\wedge\overline{dz^{n}}\otimes\frac
{d}{dz^{1}}\wedge...\wedge\frac{d}{dz^{n}}. \label{Tod31}%
\end{equation}
Since
\[
\phi_{\tau}\in C^{\infty}\left(  \text{M,}Hom\left(  \Omega_{\text{M}_{\tau}%
}^{1,0},\Omega_{\text{M}_{\tau}}^{0,1}\right)  \right)
\]
then
\[
\wedge^{n}\phi_{\tau}\in C^{\infty}\left(  \text{M,}Hom\left(  \wedge
^{n}\left(  \Omega_{\text{M}_{\tau}}^{1,0}\right)  ,\wedge^{n}\left(
\Omega_{\text{M}_{\tau}}^{0,1}\right)  \right)  \right)  \approxeq
\]%
\begin{equation}
C^{\infty}\left(  \text{M,}Hom\left(  \Omega_{\text{M}_{\tau}}^{n,0}%
,\Omega_{\text{M}_{\tau}}^{0,n}\right)  \right)  \approxeq C^{\infty}\left(
\text{M}_{\tau},\left(  \left(  \Omega_{\text{M}_{\tau}}^{n,0}\right)  ^{\ast
}\right)  ^{\otimes2}\right)  . \label{Tod32}%
\end{equation}
Combining $\left(  \ref{Tod31}\right)  $ and $\left(  \ref{Tod32}\right)  $ we
get that
\begin{equation}
\wedge^{n}\phi_{\tau}\in H^{0}\left(  \text{C,}\left(  \left(  R^{0}\pi_{\ast
}\left(  \Omega_{\mathcal{X}/\mathcal{M}\text{(M)}}^{n}\right)  \right)
^{\ast}\right)  ^{\otimes2}\right)  \label{Tod33}%
\end{equation}
is a non zero class. Moreover according to Lemma \ref{0l} $\phi_{\tau}$ is a
parallel section with respect to the pullback of the Hodge metric. This
implies that $\wedge^{n}\phi_{\tau}$ is a global parallel non zero section on
C of the line bundle
\[
\left(  \left(  R^{0}\pi_{\ast}\left(  \Omega_{\mathcal{X}/\mathcal{M}%
\text{(M)}}^{n}\right)  \right)  ^{\ast}\right)  ^{\otimes2}.
\]
From here we get that
\[
\left(  \left(  R^{0}\pi_{\ast}\left(  \Omega_{\mathcal{X}/\mathcal{M}%
\text{(M)}}^{n}\right)  \right)  ^{\ast}\right)  ^{\otimes2}%
\]
is a trivial line bundle on C. But this is impossible since we know from
\cite{LS} and the Grothendieck Riemann relative Riemann Roch Theorem that the
Chern class
\[
c_{1}\left(  \left(  R^{0}\pi_{\ast}\left(  \Omega_{\mathcal{X}/\mathcal{M}%
\text{(M)}}^{n}\right)  \right)  ^{\ast}\right)
\]
of the line bundle $\left(  R^{0}\pi_{\ast}\left(  \Omega_{\mathcal{X}%
/\mathcal{M}\text{(M)}}^{n}\right)  \right)  ^{\ast}$ is proportional to the
imaginary part of the Weil-Petersson metric. We also know that the Chern form
of the $L^{2}$ metric on $\left(  R^{0}\pi_{\ast}\left(  \Omega_{\mathcal{X}%
/\mathcal{M}\text{(M)}}^{n}\right)  \right)  ^{\ast}$ is the imaginary part of
the Weil-Petersson metric. So
\begin{equation}
\int\limits_{\text{C}}c_{1}\left(  \left(  R^{0}\pi_{\ast}\left(
\Omega_{\mathcal{X}/\mathcal{M}\text{(M)}}^{n}\right)  \right)  ^{\ast
}\right)  >0. \label{Tod34}%
\end{equation}
The inequality $\left(  \ref{Tod34}\right)  $ implies that $R^{0}\pi_{\ast
}\left(  \Omega_{\mathcal{X}/\mathcal{M}\text{(M)}}^{n}\right)  $ is not a
trivial line bundle on C. So we get a contradiction. Lemma \ref{TOD31} is
proved. $\blacksquare$

Lemma \ref{TOD31} implies Theorem \ref{TOD3}. $\blacksquare$

\textbf{Remark. }Theorem \ref{TOD3} presumably also follows from the methods
used in \cite{VZ1}, and it seems to extend to non-isotrivial morphisms
\[
\pi:\mathcal{X}_{\text{C}}\rightarrow\text{C}%
\]
whose generic fibre has a semi-ample canonical sheaf. In fact, by \cite{VZ1},
\textbf{6.2. b.}, it is sufficient to show that such a family is rigid. To
prove this, assume that there exists a deformation
\[
\pi:\mathcal{X}_{\text{C}\times\text{T}}\rightarrow\text{C}\times\text{T}%
\]
with T a projective curve, which is smooth over C$_{0}\times$T$_{0}.$ For
\[
\text{S}=\text{Y}\ominus\left(  \text{C}_{0}\times\text{T}_{0}\right)  ,
\]
according to \cite{VZ1}, \textbf{1,4,} the pullback of $\sigma^{\ast}%
S^{m}(\Omega_{\text{Y}}^{1}(\log($S$))$ under some finite map
\[
\sigma^{{\prime}}:\text{Y}^{{\prime}}\rightarrow\text{Y,}%
\]
for some $m\leq n,$ contains a big subsheaf $\mathcal{P}^{^{\prime}}.$ On the
other hand the Condition $\ref{Irr0}$ should imply that the composite
\[
\mathcal{P}^{{\prime}}\rightarrow\sigma^{\ast}S^{m}(\Omega_{\text{Y}}^{1}%
(\log(S))\rightarrow\sigma^{\ast}S^{m}(pr_{1}^{\ast}\left(  \Omega_{\text{Y}%
}^{1}(\log(S)\right)  )
\]
is non trivial. As in the proof of \cite{VZ1}, \textbf{6.5., }this would be a contradiction.

\begin{theorem}
\label{Yuk}Suppose that
\begin{equation}
\mathcal{X}\rightarrow\mathcal{C} \label{fama}%
\end{equation}
is a family of polarized CY manifold over a Riemann surfaces $\mathcal{C}$.
Suppose that the family $\left(  \ref{fama}\right)  $ contains a point
$\tau_{0}\in\mathcal{C}$ such that around it the monodromy operator
\emph{T}\ has a maximal index of unipotency, i.e.%
\begin{equation}
(\mathcal{T}^{N}-id)^{n+1}=0\text{ and }(\mathcal{T}^{N}-id)^{n}\neq0.
\label{muni}%
\end{equation}
Then the family $\left(  \ref{fama}\right)  $ is rigid.
\end{theorem}

\textbf{Proof:} We will prove that the Yukawa coupling is non-zero in an open
neighborhood of the point $\tau_{0}\in\mathcal{C}$ around which the condition
$\left(  \ref{muni}\right)  $ holds. Theorem \ref{Yuk} follows from the
existence of limit mixed Hodge structure related to the monodromy operator
$\mathcal{T}$ established in \cite{Sch}. Indeed the condition $\left(
\ref{muni}\right)  $ combined with Schmid's Theorem implies the existence of
the filtration%
\begin{equation}
\mathbb{W}_{0}\subset\mathbb{W}_{1}\subset...\subset\mathbb{W}_{2n-1}%
\subset\mathbb{W}_{2n}=\mathbb{V\otimes C}=H^{n}(\text{M}_{\tau},\mathbb{C})
\label{Filt}%
\end{equation}
such that the nilpotent operator\emph{ }%
\[
\mathcal{N}=\log(\mathcal{T}^{N}-id)
\]
acts on the filtration $\left(  \ref{Filt}\right)  $ as follows:%
\[
\mathcal{N}(\mathbb{W}_{k})\subseteq\mathbb{W}_{k-2}.
\]
In case of CY manifolds Schmid Theorem about the existence of the limit mixed
Hodge structure and condition $\left(  \ref{muni}\right)  $ imply that%
\begin{equation}
\mathcal{N}^{n}(\mathbb{W}_{2n}/\mathbb{W}_{2n-1})\approxeq\mathbb{W}_{0}.
\label{Wil}%
\end{equation}
(See \cite{Sch}.) The nilpotent orbit theorem proved in \cite{Sch} combined
with (\ref{Wil}) implies that for each point in a "small disk" containing the
point $\tau_{0}\in\mathcal{C}$ we have that
\begin{equation}
\underset{n}{\underbrace{\nabla_{\frac{\partial}{\partial\tau_{j_{1}}}}%
\circ...\circ\nabla_{\frac{\partial}{\partial\tau_{j_{n}}}}}}\omega_{\tau}%
\neq0 \label{Wil1}%
\end{equation}
and the $(0,n)$ component of
\[
\underset{n}{\underbrace{\nabla_{\frac{\partial}{\partial\tau_{j_{1}}}}%
\circ...\circ\nabla_{\frac{\partial}{\partial\tau_{j_{n}}}}}}\omega_{\tau}%
\]
is a non-zero class of cohomology. Here $\omega_{\tau}$ is a family of
holomorphic $n$ forms around the point $\tau\in\mathcal{K}$ in the Kuranishi
space that corresponds to the fibre of the point on $\mathcal{C}$ close to
$\tau_{0}.$ Now Theorem \ref{Yuk} follows directly from the definition
\ref{yuk} and Theorem \ref{TOD3}. $\blacksquare$

\subsection{Yau's Form of Schwartz Lemma and Boundedness}

S.-T. Yau proved in \cite{Y} the following theorem, which is a generalization
of the Schwarz Lemma:

\begin{theorem}
\label{Yau}Let $N$ be a complex manifold and let $h$ be a Hermitian metric on
$N$ such that the holomorphic sectional curvature $K$ of $h$ has the following
property; there exists a positive constant $c>0$ such that $K\leq-c.$ Let $C$
be either a compact or an affine Riemann surface of hyperbolic type. Let
\[
\phi:\text{C}\rightarrow N
\]
be a holomorphic map different from the constant one, then we have
\[
\phi^{\ast}(h)\leq c^{-1}g_{P},
\]
where $g_{P}$ is the Poincare metric on C, i.e. the metric with a constant curvature.
\end{theorem}

This lemma was used in \cite{Liu} to derive various height inequalities. Here
we will use the same method to prove the following Theorem:

\begin{theorem}
\label{Kef} Let
\[
\phi:\text{C}\rightarrow\overline{\mathcal{M}_{L}(\text{M})}%
\]
be a holomorphic map such that
\[
\phi(\text{C})\varsubsetneq\mathcal{D}_{\infty}%
\]
and let $x_{1},..x_{n}$ be the only n different points on C such that
\[
\phi(x_{i})\in\mathcal{D}_{\infty}.
\]
Let C$_{1}:=$C$\circleddash(x_{1}\cup...\cup x_{n}).$ Then we have the
following estimate
\[
vol(\phi(\text{C}_{1}))\leq c^{-1}|\chi(\text{C}_{1})|
\]
for the volume of $\phi(C_{1})$ in $\mathcal{M}($M$)$ with respect to the
Hodge metric $h$ on $\mathcal{M}$(M).
\end{theorem}

\textbf{Proof}: From Yau's form of the Schwarz Lemma we obtain:
\[
vol(\phi(\text{C}_{1}))=\int\limits_{C_{1}}h\leq c^{-1}\int\limits_{C_{1}%
}g_{p}.
\]
Since the Poincare metric has a constant curvature $-1$ and from the
Gauss-Bonnet theorem we obtain
\[
\int\limits_{C_{1}}g_{p}=|\chi(\text{C}_{1})|.
\]
From here Theorem \ref{Kef} follows directly. $\blacksquare$

\begin{definition}
\label{Shaf}Let $($C$;x_{1},...,x_{n})$ be a fixed Riemann surface with fix
points $x_{1},...,x_{n}.$ Let S be the divisor $x_{1}+...+x_{n}$ on C. Let us
fix a CY manifold M. We will define $Sh($C$;$S, M$)$ to be the set of all
possible families $\mathcal{X\rightarrow}$C of CY manifolds defined up to an
isomorphism over C with fixed degenerate fibres over the points $x_{1}%
,...,x_{n},$ where $\mathcal{X}$ \ is a projective manifold and the generic
fibre is a CY manifold $C^{\infty}$ equivalent to M.
\end{definition}

\begin{theorem}
\label{TOD2}Suppose that Condition $\ref{Irr0}$ or its equivalent $\left(
\ref{Be1}\right)  $ holds for some fibre M$_{\tau}$ of a family of CY
manifolds over a fix Riemann surface C with fixed points if degenerations.
Then the set $Sh($C$;$S$)$ is finite.
\end{theorem}

\textbf{Proof}: The proof of Theorem \ref{TOD2} will be done in two steps. The
first step will be to prove that the set of maps $Sh($C$;$S$)$ is discrete.
The second step is to show that for each $\phi\in Sh($C$;$S$)$ the volume of
$\phi($C$)$ is bounded by a universal constant. Then Theorem \ref{TOD2} will
follow directly from the theorem of Bishop proved in \cite{Bi}, which implies
that the set $Sh$(C,S) is compact.

\begin{lemma}
\label{TOD21}The set of maps $Sh$(C,S) is discrete.
\end{lemma}

\textbf{Proof: }Since we assumed that the CY manifold M satisfies Condition
$\left(  \ref{Irr}\right)  $, Lemma \ref{TOD21} follows directly from Theorem
\ref{TOD3}. $\blacksquare$

We will need the following Lemma for establishing the second step described above:

\begin{lemma}
\label{t1} Let $Sh$(C,S) be the set defined in Definition \ref{Sh}, where S is
the divisor $x_{1}+...+x_{n}$ on C. Then $Sh$(C,S) is compact.
\end{lemma}

\textbf{Proof}: Each point $z\in Sh$(C,S) is represented by a map
\[
\phi_{z}:\text{C}\rightarrow\mathcal{T}(\text{M})/\Gamma
\]
which satisfies the conditions of Definition \ref{Sh}. Let $g$ be the
K\"{a}hler (1,1) form of the Hodge metric on $\mathcal{T}($M$)/\Gamma$, then
we can define a height function $h$ on $Sh$(C,S) as follows:%

\begin{equation}
h(\phi_{z})=\int\limits_{\phi(\text{C}_{z})}\operatorname{Im}g=vol(\phi
_{z}(\text{C})). \label{f1}%
\end{equation}
Since the metric $g$ is a K\"{a}hler metric with logarithmic growth on X, then%

\begin{equation}
h(\phi_{z})=\int\limits_{\overline{\mathcal{T}(\text{M})/\Gamma}%
}\operatorname{Im}\left(  g\right)  \wedge\mathcal{P}(\phi(\text{C))}%
=vol(\phi_{z}(\text{C})). \label{f2}%
\end{equation}
where $\mathcal{P}(\phi_{z}($C$))$ is the Poincare dual of the homology class
of $\overline{\phi_{z}(\text{C})}$ in $\overline{\mathcal{T}(\text{M})/\Gamma
}.$ Since $\phi_{z}\in Sh$(C,S) then it is a deformation of the map f. So for
all $\phi_{z}\in Sh$(C,S)$,$ $\overline{\phi_{z}(\text{C})}$ realizes a fixed
class of homology in $\overline{\mathcal{T}(\text{M})/\Gamma}$. The integrals
that appeared in $\left(  \ref{f1}\right)  $ and $\left(  \ref{f1}\right)  $
are finite. This follows from Theorem 5.1. proved in \cite{LS}. Theorem
\ref{Kef} implies that $h(\phi_{z})$ is a bounded function on $Sh$(C,S).

We will prove now that $Sh$(C,S) is a compact set. Since $Sh$(C,S) is a
discrete set, the compactness of $Sh$(C,S) will imply that it is a finite set.
The compactness of $Sh$(C,S) will follows if we prove that from any sequence
$\{\phi_{n}\}$ in $Sh$(C,S) there exists a subsequence which converges weakly
to an algebraic map, i.e. the corresponding subsequence of images of S
converges to an algebraic subvariety of X. Bishop's theorem implies that. See
page 292 of \cite{Bi} and also the Appendix II. Lemma \ref{t1} is proved.
$\blacksquare$

Theorem \ref{TOD2} follows directly from Lemmas \ref{TOD21} and \ref{t1}.
$\blacksquare$

\begin{theorem}
\label{Yau1}Suppose that M is a CY manifold such that each non isotrivial
family of CY manifolds over a fixed Riemann surface with fixed points of
degenerations and each fibre appears as a fibre in a fixed component of the
Teichm\"{u}ller space of M is rigid. Then Shafarevich conjecture holds for
that type of CY manifolds.
\end{theorem}

\section{Higgs Bundles, VHS and Rigidity}

\subsection{Non-Rigid Families of CY Manifolds and Monodromy}

We are going to study the relations between Higgs bundles, Variations of Hodge
Structures (VHS) and the existence of non rigid families of CY manifolds.

Let S$^{0}$ and T$^{0}$ be smooth quasi-projective manifolds. Let
\begin{equation}
\pi:\mathcal{X\rightarrow}\text{S}^{0}\times\text{T}^{0}=\mathcal{Y}_{0}
\label{Den}%
\end{equation}
be a non rigid family of Calabi-Yau 3-folds, such that the induced map into
the moduli space is generically finite. Let
\[
\mathbb{V}=R^{3}\pi_{\ast}\mathbb{Z}%
\]
denote the flat bundle associated with the the third cohomology group. The
Variation of Hodge Structures (VHS) of weight 3 of $\pi$ associated with the
family $\left(  \ref{Den}\right)  $ will be $\mathbb{V\otimes C}.$ The
corresponding Hodge bundles are then
\[
\mathbb{V\otimes C}=\pi_{\ast}\Omega_{\mathcal{X}\text{/}\mathcal{Y}_{0}}%
^{3}\text{ }\oplus R^{1}\pi_{\ast}\Omega_{\mathcal{X}\text{/}\mathcal{Y}_{0}%
}^{2}\oplus R^{2}\pi_{\ast}\Omega_{\mathcal{X}\text{/}\mathcal{Y}_{0}}%
^{1}\oplus R^{3}\pi_{\ast}\mathcal{O}_{\mathcal{X}}.
\]
The cup product of the Kodaira-Spencer class of $\pi$%
\[
\theta_{0}^{p,q}:\ E_{0}^{p,q}\rightarrow E_{0}^{p,q}\otimes\Omega
_{\mathcal{Y}_{0}}^{1},
\]
where
\[
E_{0}^{3,0}=\pi_{\ast}\Omega_{\mathcal{X}\text{/}\mathcal{Y}_{0}}^{3},\text{
}E_{0}^{2,1}=R^{1}\pi_{\ast}\Omega_{\mathcal{X}\text{/}\mathcal{Y}_{0}}%
^{2},\text{ }E_{0}^{1,2}=R^{2}\pi_{\ast}\Omega_{\mathcal{X}\text{/}%
\mathcal{Y}_{0}}^{2}\text{ and }E_{0}^{0,3}=R^{3}\pi_{\ast}\mathcal{O}%
_{\mathcal{X}}.
\]
will provide
\[
\left(  \underset{p+q=3}{\oplus}E_{0}^{p,q},\underset{p+q=3}{\oplus}\theta
_{0}^{p,q}\right)
\]
with the Higgs bundle structure on $\mathbb{V}$ corresponding to the VHS$.$

\begin{theorem}
\label{Zuo2}Suppose that S$^{0}$ and T$^{0}$ are Zariski open sets in
projective varieties S and T. Suppose that
\begin{equation}
\pi:\ \mathcal{X\rightarrow Y}_{0}=\text{S}^{0}\times\text{T}^{0} \label{CYY}%
\end{equation}
is a non rigid family of Calabi-Yau manifolds, such that the induced map to
the moduli space is generically finite. Let
\[
\mathbb{V\rightarrow}\mathcal{Y}_{0}=\text{S}^{0}\times\text{T}^{0}%
\]
be the VHS induced by the family $\left(  \ref{CYY}\right)  $. Suppose that
$(s_{i},t_{j})$ is a point of
\[
\mathfrak{D}_{\infty}=\left(  \text{S}\ominus\text{S}^{0}\right)
\times\left(  \text{T}\ominus\text{T}^{0}\right)
\]
Let $\gamma$\textsl{$_{s_{i}}$ }and $\gamma$\textsl{$_{t_{j}}$ }be two short
loops around the divisors T$_{s_{i}}$ and S$_{t_{j}}$. Let S$_{t_{0}}$ be the
curve in S-direction passing through the point (s$_{0},$t$_{0}$). \textsl{
}Suppose that the monodromy operators $\rho(\gamma_{s_{i}})=\mathcal{T}_{i}$
and $\rho(\gamma_{t_{j}})=\mathcal{T}_{j}$ have infinite order, then the
endomorphisms%
\[
\mathcal{T}_{i}:\ \mathbb{V}|_{(s_{0},t_{0})}\rightarrow\mathbb{V}%
|_{(s_{0},t_{0})}%
\]
and%
\[
\mathcal{T}_{j}:\ \mathbb{V}|_{(s_{0},t_{0})}\rightarrow\mathbb{V}%
|_{(s_{0},t_{0})}%
\]
are distinct.
\end{theorem}

\textbf{Proof: }The idea of the proof is to construct on the restrictions of
VHS $\mathbb{V}|_{\text{S}^{0}}$ and $\mathbb{V}|_{\text{T}^{0}}$ two
endomorphisms by using the fact that the family of CY manifold $\mathcal{X}%
\rightarrow$S$^{0}\times$T$^{0}$ is not rigid. According to Lemma \ref{0l}
these two endomorphisms are parallel sections of the Hodge bundle $R^{1}%
\pi_{\ast}\Omega_{\mathcal{X}/\text{S}^{0}\times\text{T}^{0}}^{n-1}.$ Theorem
\ref{Zuo2} will follow from the way the weight filtrations of nilpotent
endomorphisms introduced by W. Schmid are constructed, the fact that the two
filtrations will be different and the Theorem about the complete reducibility
of local systems of VHS over quasi-projective manifolds due to
Deligne.\newline

There are two ways to obtain endomorphisms on $\mathbb{V}|_{\text{S}_{t}}$ and
$\mathbb{V}|_{\text{T}_{j}},$ where S$_{i}:=$S$\times t_{0}$ and T$_{j}%
=s_{0}\times$T$.$\newline

\textbf{1. }Consider the two projections%
\[
p_{\text{S}}:\text{S}\times\text{T}\rightarrow\text{S and }p_{\text{T}%
}:\text{S}\times\text{T}\rightarrow\text{T}.
\]
One has%

\[
\Omega_{\text{S}\times\text{T}}^{1}(\log\mathcal{D}_{\text{S}\times\text{T}%
})=p_{\text{S}}^{\ast}\left(  \Omega_{\text{S}}^{1}(\log\mathcal{D}_{\text{S}%
})\right)  \oplus p_{\text{T}}^{\ast}\left(  \Omega_{\text{T}}^{1}%
(\log\mathcal{D}_{\text{T}})\right)
\]
where $\mathcal{D}_{\text{S}\times\text{T}}$ is the discriminant locus of the
family $\left(  \ref{CYY}\right)  .$ Then the restriction of the Higgs map to
S$_{t}$ defines a natural map%

\begin{equation}
\theta|_{\text{S}_{t}}:(p_{\text{S}}^{\ast}\left(  \Theta_{\text{S}}%
(-\log\mathcal{D}_{\text{S}}\right)  )\oplus p_{\text{T}}^{\ast}\left(
\Theta_{\text{T}}(-\log\mathcal{D}_{\text{T}})\right)  _{\text{S}_{t}%
}\rightarrow End(\mathbb{V})|_{\text{S}_{t}}, \label{kang}%
\end{equation}
\newline by identifying the sections of $(p_{\text{S}}^{\ast}\left(
\Theta_{\text{S}}(-\log\mathcal{D}_{\text{S}}\right)  )\oplus p_{\text{T}%
}^{\ast}\left(  \Theta_{\text{T}}(-\log\mathcal{D}_{\text{T}})\right)
_{\text{S}_{t}}$ with Kodaira-Spencer classes. Note that
\[
p_{\text{T}}^{\ast}\Theta_{\text{T}}(-\log\mathcal{D}_{\text{T}}%
)|_{\text{S}_{t}}\simeq\mathcal{O}_{\text{S}_{t}}\text{ and }p_{\text{S}%
}^{\ast}\Theta_{\text{S}}(-\log\mathcal{D}_{\text{S}})|_{\text{T}_{j}}%
\simeq\mathcal{O}_{\text{T}_{j}}.
\]
Let
\[
1_{T}\in p_{\text{T}}^{\ast}\Theta_{\text{T}}(-\log\mathcal{D}_{\text{T}%
})|_{\text{S}_{i}}\text{ and }1_{S}\in p_{\text{S}}^{\ast}\Theta_{\text{S}%
}(-\log\mathcal{D}_{\text{S}})|_{\text{T}_{j}}%
\]
be the constant sections. Then by using a Theorem of Jost and Yau proved in
\cite{JY} one obtains endomorphisms%

\begin{equation}
\theta|_{S_{t}}(1_{\text{S}t}):\ \mathbb{V}_{\text{S}_{i}^{0}}\rightarrow
\mathbb{V}_{\text{S}_{i}^{0}}\text{ and }\theta|_{\text{T}_{j}}(1_{\text{T}%
_{t}}):\mathbb{V}_{\text{T}_{j}^{0}}\rightarrow\mathbb{V}_{\text{T}_{j}^{0}}.
\label{kang3}%
\end{equation}
See also \cite{Zuo} for another proof. This proof is based on the observation
that the image of the map $\left(  \ref{kang}\right)  $ is contained in the
kernel of the induced Higgs map on $\text{End}(\mathbb{V}).$ Therefore, the
image is a Higgs subsheaf with the trivial Higgs field. The Higgs
poly-stability of $\text{End}(\mathbb{V})|_{S_{t}}$ implies that any section
in this subsheaf is flat. One sees also that this flat section is of Hodge
type (-1,1).\newline

\textbf{2.} Let $\mathbb{V}$ be an arbitrary Variation of Hodge Structure over
some quasi-projective space by taking the local monodromy by the following
consideration. Let S$^{0}=\,$S$\setminus\{s_{0},,,s_{m}\}$ and T$^{0}%
=$T$\setminus\{t_{0},...,t_{n}\}$ be two smooth projective curves. The local
system $\mathbb{V}$ corresponds to a representation%

\[
\rho:\pi_{1}(\text{S}^{0},\ast)\times\pi_{1}(\text{T}^{0},\ast)\rightarrow
GL(V),
\]
and the restriction $\rho|_{\text{S}_{t}}$ corresponds to the restricted representation%

\[
\rho|_{\pi_{1}(\text{S}^{0},\ast)}:\ \pi_{1}(\text{S}^{0},\ast)\rightarrow
GL(V).
\]
\newline Let $\gamma_{t_{j}}$ be a short loop around some point $t_{j},$ then
it commutes with $\pi_{1}($S$^{0},\ast).$ Hence, $\rho(\gamma_{t})$ descends
to an endomorphism, say%

\[
\rho(\gamma_{t_{j}}):\ \mathbb{V}_{\text{S}_{t}^{0}}\rightarrow\mathbb{V}%
_{\text{S}_{t}^{0}}.
\]
\newline Similarly, one gets also an endomorphism $\rho(\gamma_{s_{i}})$ on
$\mathbb{V}_{\text{T}_{s}}\rightarrow\mathbb{V}_{\text{T}_{s}}.$ It is well
known that those two endomorphisms $\rho(\gamma_{t_{j}})$ and $\rho
(\gamma_{s_{i}})$ are unipotent. Let
\[
\mathcal{N}_{s_{i}}=\log(\mathcal{T}_{i}-id)\text{ }\text{and}\text{
}\mathcal{N}_{t_{j}}=\log(\mathcal{T}_{j}-id)
\]
be the nilpotent part of $\mathcal{T}_{i}$ and $\mathcal{T}_{j}$ respectively
which are non-trivial by assumption. Let%
\[
\theta|_{S_{t}}(1_{\text{S}_{t}}):\ \mathbb{V}_{\text{S}_{i}^{0}}%
\rightarrow\mathbb{V}_{\text{S}_{i}^{0}}%
\]
be the endomorphism of the variations of the Hodge structure induced by
$\theta|_{S_{t}}(1_{\text{S}t})$ in $\left(  \ref{kang3}\right)  $. We noticed
that $\theta|_{S_{t}}(1_{\text{S}t})$ and $\theta|_{T_{j}}(1_{\text{T}s})$ are
nilpotent. Similar to what Schmid did for nilpotent endomorphisms on vector
space we see that $\theta|_{S_{t}}(1_{\text{S}t})$ induces a weight filtration
of the local system%

\[
0\subset\mathbb{W}_{0}\subset\mathbb{W}_{1}\subset...\mathbb{W}_{2k-1}%
\subset\mathbb{W}_{2k}=\mathbb{V}_{\text{S}_{t_{0}}^{0}},
\]
such that
\[
\theta|_{S_{t}}(1_{\text{S}t})(\mathbb{W}_{l})\subset\mathbb{W}_{l-2}.\newline%
\]
Since the local system $\mathbb{V}_{S_{t_{0}}^{0}}$ is completely reducible by
Deligne, one gets an isomorphism
\[
\mathbb{V}|_{s_{j}}=\underset{j}{\oplus}\mathbb{W}_{j}/\mathbb{W}_{j-1}.
\]
The isomorphism $\theta|_{S_{t}}(1_{\text{S}t})$ shifts $\mathbb{W}%
_{l}/\mathbb{W}_{l-1}$ to $\bigoplus_{j\leq l-2}\mathbb{W}_{j}/\mathbb{W}%
_{j-1}.$ On the other hand, regarding $\mathcal{T}_{i}$ as an element in the image%

\[
\rho:\ \pi_{1}(\text{S}_{t_{0}},s_{0})\rightarrow GL(V)
\]
the above decomposition on the fibre $\mathbb{V}_{s_{0},t_{0}}$ decomposes
$\mathcal{N}_{s_{i}}=\log(\mathcal{T}_{i}-id)$ as a direct sum of
endomorphisms on $\left(  \mathbb{W}_{j}\right)  _{(s_{0},t_{0})}/\left(
\mathbb{W}_{j-1}\right)  _{(s_{0},t_{0})}.$ In particular, $\mathcal{N}%
_{s_{i}}$ can not be equal to $\mathcal{N}_{t_{j}}$ as endomorphism on the
fibre $\mathbb{V}_{(s_{0},t_{0})}$ since the filtrations of those two parallel
endomorphisms of VHS will be different since the families obtained from the
restrictions of the family $\left(  \ref{CYY}\right)  $ on $S$ and $T$ are not
isomorphic. This follows from local Torelli Theorem. The two filtrations will
be invariant under the actions of $\log(\mathcal{T}_{i}-id)$ and
$\log(\mathcal{T}_{j}-id)$ respectively. Theorem \ref{Zuo2} is proved.
$\blacksquare$

\begin{corollary}
\label{Zuo1}Suppose that S$^{0}$ and T$^{0}$ are Zariski open sets in
projective varieties S and T. Suppose that
\[
\pi:\mathcal{X\rightarrow Y}=S\times T
\]
is the family of Calabi-Yau manifolds $\left(  \ref{CYY}\right)  $, such that
the induced map to the moduli space $\overline{\mathcal{M}_{L}(\text{M})}$ is
generically finite. Suppose that the discriminant locus $\mathfrak{D}_{\infty
}=$ $\overline{\mathcal{M}_{L}(\text{M})}\ominus$ $\mathcal{M}_{L}($M$)$ is a
irreducible divisor and at some point $(s_{0},t_{0})$ of
\[
\mathfrak{D}_{\infty}=\left(  \text{S}\ominus\text{S}^{0}\right)
\times\left(  \text{T}\ominus\text{T}^{0}\right)
\]
the local monodromy operators of the restrictions of the family $\left(
\ref{CYY}\right)  $ on S and T around the points $s_{0}\in$S and $t_{0}\in$T
are infinite$.$ Then the family $\left(  \ref{CYY}\right)  $ is rigid.
\end{corollary}

\textbf{Proof: }Since we assumed that the discriminant locus
\[
\mathfrak{D}_{\infty}=\overline{\mathcal{M}_{L}(\text{M})}\ominus
\mathcal{M}_{L}(\text{M})
\]
is an irreducible divisor in $\overline{\mathcal{M}_{L}(\text{M})}$ then the
monodromy operators acting on the middle cohomology of a fixed CY manifold
induced by any loops $\gamma$ in $\mathcal{M}_{L}($M$)$ around the image of
the point $(s_{0},t_{0})$ in $\mathfrak{D}_{\infty}$ will be the same. Since
we assumed that the monodromy operators are infinite we get a contradiction
with Theorem\ \ref{Zuo2}. $\blacksquare$

\begin{remark}
It is a well known fact that the discriminant locus $\mathfrak{D}_{\infty}$ in
the Hilbert scheme of hypersurfaces in $\mathbb{CP}^{n}$ for $n\geq2$ is an
irreducible divisor.
\end{remark}

The motivation of the next Theorem is based on a question raised by S.-T. Yau
and this question is closely related to the example constructed in
\textbf{Appendix I. }We will use the following definitions in the formulation
of Theorem \ref{Zuo}: We will say that $\mathbb{V}$ is a K3 like, abelian or
elliptic local system it is a VHS of weight two induced by an Euclidean
lattice with signature $(2,k)$ where $0<k\leq19,$ or a VHS of an abelian
variety or elliptic curve respectively.

\begin{theorem}
\label{Zuo}Suppose that
\begin{equation}
\pi:\ \mathcal{X\rightarrow}\text{S}^{0}\times\text{T}^{0}=\mathcal{Y}_{0}
\label{nrf}%
\end{equation}
is a family of Calabi-Yau manifolds, such that the induced map into the moduli
space is generically finite where S$^{0}$ and T$^{0}$ are Zariski open sets in
projective varieties S and T. Suppose that
\[
\mathbb{V\rightarrow}\text{S}^{0}\times\text{T}^{0}=\mathcal{Y}_{0}%
\]
is the VHS induced by the family $\left(  \ref{CYY}\right)  $. Fixing t$\in
$T$^{0}$, then the local system $\mathbb{V}\otimes\mathbb{C}|_{\text{S}_{t}}$
has a splitting over $\mathbb{C}$ and the splittings of the local systems are
either
\begin{equation}
\mathbb{V}\otimes\mathbb{C}|_{\text{S}_{t}}=\mathbb{G}_{0}\oplus\mathbb{G}%
_{1}\oplus\mathbb{G}_{2} \label{K31}%
\end{equation}
where $\mathbb{G}_{0}$ and $\mathbb{G}_{2}$ are K3 liked local systems which
are isomorphic to each other and $\mathbb{G}_{1}$ is unitary local system or
\begin{equation}
\mathbb{V}\otimes\mathbb{C}|_{\text{S}_{t}}=\mathbb{G}_{0}\oplus\mathbb{G}%
_{1}\oplus\mathbb{G}_{2}\oplus\mathbb{G}_{3}\oplus\mathbb{G}_{4} \label{Ell1}%
\end{equation}
where $\mathbb{G}_{0}$ and $\mathbb{G}_{4}$ are isomorphic elliptic-like local
systems, $\mathbb{G}_{2}$ is an abelian-like local system and $\mathbb{G}_{1}$
and $\mathbb{G}_{3}$ are unitary local systems that are isomorphic to each other.
\end{theorem}

\textbf{Proof:} Consider the endomorphism%

\[
\sigma:\ \mathbb{V}\otimes\mathbb{C}|_{S_{t}^{0}}\rightarrow\mathbb{V}%
\otimes\mathbb{C}|_{S_{t}^{0}},
\]
which is the image of%

\[
1_{S_{t}}:\mathcal{O}_{\text{S}_{t}}\rightarrow p_{\text{T}}^{\ast}%
\Theta_{\text{T}}(-\log\mathcal{D}_{\text{T}})|_{\text{S}_{t}}\rightarrow
\text{End}(\mathbb{V}|_{S_{t}}).
\]
Since we assumed that the family $\left(  \ref{nrf}\right)  $ is non-rigid
then we have that $\sigma^{3}=0.$ Otherwise if $\sigma^{3}\neq0$ then by
Theorem \ref{Yuk} the family $\left(  \ref{nrf}\right)  $\ will be rigid.
Clearly $\sigma$ is of (-1,1) type. We know from Lemma \ref{0l} that any non
rigid deformation corresponds to a parallel section of the pullback of the
tangent bundle of the moduli space. Thus this parallel section can be
identified with some Kodaira-Spencer class $\theta$. Applying Cor. 6.5 in
\cite{VZ1} to non-rigid families of CY 3-folds, the three times iterated
Kodaira-Spencer map%

\[
(i_{T_{s}}^{\ast}\theta)^{3}:\ \mathbb{V}|_{T_{s}}\rightarrow\mathbb{V}%
_{T_{s}}\otimes S^{3}\Omega_{T_{s}}^{1}(\log D_{T_{s}}),\quad\forall s\in
S^{0}%
\]
must be zero. Fixing $t\in T,$ and varying $s$, the restriction $(i_{T_{s}%
}^{\ast}\theta)|_{t,s}$ is precisely $\sigma.$ Next we will consider two
different cases:\newline\textbf{Case 1.} Suppose that $\sigma^{2}=0.$ We
consider the weight filtration studied in \cite{Sch} defined by $\sigma$.%

\[
\mathbb{W}_{0}\subset\mathbb{W}_{1}\subset\mathbb{W}_{2}=\mathbb{V}|_{S_{t}}.
\]
We set%

\[
\mathbb{W}_{0}=\sigma\left(  \mathbb{W}_{2}\right)  ,\quad\mathbb{W}%
_{1}=\text{ker}\left(  \sigma:\mathbb{W}_{2}\rightarrow\mathbb{W}_{0}\right)
.
\]
One has $\sigma:\mathbb{W}_{2}/\mathbb{W}_{1}\simeq\mathbb{W}_{0}.$ Note that
the Hodge bundle corresponding to $\mathbb{W}_{0}$ has the form%
\[
\mathbb{W}_{0}=\sigma\left(  E^{3,0}\right)  \oplus\sigma\left(
E^{2,1}\right)  \oplus\sigma\left(  E^{1,2}\right)  ,
\]
setting%

\[
\mathbb{G}_{0}=\mathbb{W}_{0},\text{ }\mathbb{G}_{2}=\mathbb{W}_{2}%
/\mathbb{W}_{1}%
\]
we obtain
\[
\sigma:\ \mathbb{G}_{2}\simeq\mathbb{W}_{0}=\mathbb{G}_{0}.
\]
Since from the definition of the mixed Hodge structure it follows that the
pure Hodge structure on $\mathbb{W}_{2}/\mathbb{W}_{1}$ is of weight two then
the VHS on $\mathbb{G}_{0}$ and $\mathbb{G}_{2}$ are K3 liked local system
since the fibre of the family $\left(  \ref{nrf}\right)  $ is a CY manifold
and so%
\[
\dim_{\mathbb{C}}\sigma\left(  E^{3,0}\right)  =1.
\]
\newline Note that the Hodge bundle corresponding to $\mathbb{W}_{1}$ has the form%

\[
\text{ker}(\sigma:E^{2,1}\rightarrow E^{1,2})\oplus\text{ker}(\sigma
:E^{1,2}\rightarrow E^{0,3})\oplus E^{0,3}.
\]
The Hodge bundle corresponding to $\mathbb{G}_{1}=\mathbb{W}_{1}%
/\mathbb{W}_{0}$ has the form%

\[
\text{ker}(\sigma:E^{2,1}\rightarrow E^{1,2})/\sigma\left(  E^{3,0}\right)
\oplus\text{ker}(\sigma:E^{1,2}\rightarrow E^{0,3})/\sigma\left(
E^{2,1}\right)  .
\]
Since $\mathbb{V}|_{S_{t}}$ is complete reducible by Deligne, we obtain%

\[
\mathbb{V}|_{S_{t}}\simeq\bigoplus_{i=0}^{2}{}\mathbb{G}_{i}.
\]
\newline\textbf{Case 2)} Suppose now $\sigma^{2}\not =0,\,\sigma^{3}=0.$
Consider again the weight filtration%

\[
0\subset\mathbb{W}_{0}\subset\mathbb{W}_{1}\subset\mathbb{W}_{2}%
\subset\mathbb{W}_{3}\subset\mathbb{W}_{4}=\mathbb{V}|_{S_{t}}.
\]
\newline Set%

\[
\mathbb{W}_{0}:=\sigma^{2}(\mathbb{W}_{4}),\quad\mathbb{W}_{3}:=\text{ker}%
\left(  \sigma^{2}:\mathbb{W}_{4}\rightarrow\mathbb{W}_{0}\right)  .
\]
Note that $\sigma$ is of (-1,1)-type, so the only possible non-zero parts of
$\sigma^{2}$ are%

\[
\sigma^{2}(E^{2,1})\subseteq E^{0,3},\quad\sigma^{2}\left(  E^{3,0}\right)
\subset E^{1,2}%
\]
\newline they are dual to each other. Thus, if $\sigma^{2}\not =0,$ then%

\[
\sigma^{2}(E^{2,1})=E^{0,3}\text{ and}\quad\sigma^{2}(E^{3,0})\subset
E^{1,2}\,.
\]
Therefore,%
\[
\sigma^{2}:\ \mathbb{W}_{4}/\mathbb{W}_{3}\simeq\mathbb{W}_{0}.
\]
And the Hodge bundle corresponding to $\mathbb{W}_{0}$ has the form%

\[
\theta:\ F_{0}^{1,0}\rightarrow F_{0}^{0,1},
\]
\newline where
\[
F_{0}^{1,0}=\sigma^{2}(E^{3,0}),\,F_{0}^{0,1}=E^{0,3}%
\]
and $\theta$ is the original Kodaira-Spencer class restricted to%

\[
\bigoplus_{p+q=1}F_{0}^{p,q}\subset\bigoplus_{p+q=3}E^{p,q}.
\]
\newline Thus the VHS on $\mathbb{W}_{0}$ is like the VHS on an elliptic curve.

We set now%

\[
\mathbb{W}_{3}=\text{ker}\left(  \sigma^{2}:\mathbb{W}_{4}\rightarrow
\mathbb{W}_{0}\right)
\]
\newline whose Hodge bundle has the form%

\[
\mathbb{W}_{3}=\text{ker}(\sigma^{2}:\ E^{2,1}\rightarrow E^{0,3})\oplus
E^{1,2}\oplus E^{0,3}.
\]
Let $\mathbb{W}_{1}=\sigma(\mathbb{W}_{3}).$ Then $\mathbb{W}_{1}$ as a Hodge
bundle has the form%

\[
\mathbb{W}_{1}=\sigma(\text{ker}(\sigma^{2}:E^{2,1}\rightarrow E^{0,3}%
))\oplus\sigma(E^{1,2})=\sigma(\text{ker}(\sigma^{2}:E^{2,1}\rightarrow
E^{0,3}))\oplus E^{0,3}.
\]
Let $\mathbb{W}_{2}=\{v\in|\sigma(v)\in\mathbb{W}_{0}\}.$ Clearly we have%

\[
\sigma:\mathbb{W}_{3}/\mathbb{W}_{2}\simeq\mathbb{W}_{1}/\mathbb{W}_{0}.
\]
Remember that the Hodge bundle corresponding to $\mathbb{W}_{0}$ has weight
one and thus the form%

\[
\sigma^{2}(E^{3,0})\oplus E^{0,3}.
\]
The Hodge bundle corresponding to $\mathbb{W}_{1}/\mathbb{W}_{0}$ has the form%

\[
F^{0,0}=\sigma(\text{ker}(\sigma^{2}:E^{2,1}\rightarrow E^{0,3})/\sigma
^{2}(E^{3,0}).
\]
So, we see that%

\[
\sigma:\mathbb{W}_{3}/\mathbb{W}_{2}\simeq\mathbb{W}_{1}/\mathbb{W}_{0}%
\]
are unitary local system.

Finally we want to determine the Hodge bundle corresponding to $\mathbb{W}%
_{2}$. It is just the praimage of the Hodge bundle corresponding to
$\mathbb{W}_{0}$, namely%

\[
\sigma^{-1}(\sigma^{2}(E^{3,0}))\cap E^{2,1}\oplus E^{1,2}\oplus E^{0,3}.
\]
So, the Hodge bundle corresponding to $\mathbb{W}_{2}/\mathbb{W}_{1}$ has the form%

\[
\sigma^{-1}(\sigma^{2}(E^{3,0}))\cap E^{2,1}\oplus E^{1,2}/\sigma
(\text{ker}(\sigma^{2}:E^{2,1}\rightarrow E^{0,3})).
\]
It is easy to see that its is again abelian variety liked local system.
Theorem \ref{Zuo} is proved. $\blacksquare$\

\section{\label{cexp} Appendix I. Complex Exponential Maps in K\"{a}hler
Geometry}

It is stated in \cite{GH} that the following conditions are equivalent:

\begin{enumerate}
\item \textbf{\ }$g_{i,\overline{j}}$ is a K\"{a}hler metric on a complex
manifold M.

\item The $(1,1)$ form $\omega:=\operatorname{Im}(g)$ is closed.

\item The complex structure operator $J_{M}$ on M is a parallel tensor with
respect to the Levi-Cevita connection, i.e. $\nabla J_{\text{M}}=0.$

\item Around any point $m_{0}\in M$ there exist holomorphic coordinates
$\{z^{1},...,z^{n}\}$ in an open set $\mathcal{U\subset}M$ such that locally
the metric $g$ is given by:
\[
g_{i,\overline{j}}=\delta_{ij}+\frac{1}{4}\sum_{k,l}R_{i\overline
{j},k\overline{l}}z^{k}\overline{z^{l}}+...
\]
where $R_{i\overline{j},k\overline{l}}$ is the curvature tensor.
\end{enumerate}

The coordinates $\{z^{1},...,z^{n}\}$ will be called flat coordinates with
respect to the K\"{a}hler metric g.

\subsection{A Geometric Construction of the Flat Coordinates of a K\"{a}hler
Metric}

Let us fix a point $m_{0}$ in the complex K\"{a}hler manifold M. Let
$T_{m_{0}}$ be the tangent space at the point $m_{0}\in M.$ We will consider
for the moment the tangent space $T_{m_{0}}$ as a real 2n dimensional vector
space. Let $e_{1}\in$ $T_{m_{0}}$ be a vector of length 1. We will define
$e_{n+1}:=J_{\text{M}}e_{1}.$ Let $e_{2}$ be a vector perpendicular to the
vectors $e_{1}$ and \ $e_{n+1}$ and $\left\Vert e_{2}\right\Vert =1.$ We will
define $e_{n+2}:=J_{\text{M}}e_{2}.$ Continuing this process we obtain a basis
in $T_{m_{0}}$ consisting of vectors
\[
\{e_{1},..,e_{n},J_{\text{M}}e_{1}=e_{n+1},...,J_{\text{M}}e_{n}=e_{2n}\}
\]
such that they satisfy $\left\langle e_{i},e_{j}\right\rangle =\delta_{ij}%
.\,$Let $\gamma_{i}(t)$ be geodesics for $|t|<\varepsilon$ with respect to the
metric g on M such that
\[
\gamma_{i}(0)=m_{0}\in\text{M}%
\]
and
\[
\frac{d\gamma_{i}(t)}{dt}\left\vert _{t=0}\right.  =e_{i}%
\]
\textit{for} $i=1,..,n.$

\begin{definition}
\label{distr}We will define two dimensional distributions $\mathcal{D}_{i}$ in
the tangent bundle $T($M$)$ in a small neighborhood $\mathcal{U}$ of the point
$m_{0}\in$M$.$ $\mathcal{D}_{i}$ by parallel transportation along the
geodesics containing the point $m_{0}$ of the two dimensional subspaces
$\mathcal{E}_{i}(0)\subset T_{m_{0},\text{M}}$ span by $e_{i}$ and
$J_{\text{M}}e_{i}$ for $i=1,...,n.$
\end{definition}

\begin{theorem}
\label{g2} Locally around the point $m_{0}\in$M there exist one dimensional
complex manifolds $Z_{i}$ for each $1\leq i\leq n$ such that
\end{theorem}

\begin{enumerate}
\item \textit{For each point }$x\in$\textit{\ }$\gamma_{i}(t)\subset Z_{i}$
\textit{\ the tangent space }$T_{x,Z_{i}}$ \textit{of }$Z_{i}$\textit{\ is
equal to }$\mathcal{D}_{i}(t)$\textit{. }

\item $Z_{i}$\textit{\ is a totally geodesic two dimensional real
submanifold.}
\end{enumerate}

\textbf{Proof: }We will construct $Z_{i}$ for $i=1,...,n$ based on the
following Proposition:

\begin{proposition}
\label{g32}The distributions $\mathcal{D}_{i}$ are integrable distributions.
\end{proposition}

\textbf{Proof:} Let us denote by $\widetilde{e}_{i}$ $\ $and $\widetilde
{e}_{i+n}$ the parallel transported vector fields
\[
e_{i}\text{ }\&\text{ }e_{i+n}\in T_{m_{0},\text{M}}%
\]
along the geodesics passing through $m_{0}\in$M for $i=1,...,n.$ Clearly,
$\widetilde{e}_{i}$\& $\widetilde{e}_{i+n}\in\mathcal{D}_{i}.$ The
integrability of \ the distribution $\mathcal{D}_{i}$ in the tangent bundle
$T($M$)$ will follow if we show that the Lie bracket $[\widetilde{e_{i}%
},\widetilde{e_{i+n}}]$ fulfills the condition $[\widetilde{e}_{i}%
,\widetilde{e}_{i+n}]\in\mathcal{D}_{i}.$ At the point $m_{0}$ we have
\begin{equation}
\lbrack\widetilde{e}_{i},\widetilde{e}_{i+n}]|_{m_{0}}=0\in\mathcal{D}_{i}.
\label{0}%
\end{equation}
The fact that the metric $g$ is a K\"{a}hler implies that $\nabla J_{\text{M}%
}=0.$ Using this fact and that $e_{i+n}=J_{\text{M}}e_{i},$ we deduce the
following formula:
\begin{equation}
\nabla\lbrack\widetilde{e}_{i},\widetilde{e}_{i+n}]=\nabla\lbrack\widetilde
{e}_{i},J_{\text{M}}\widetilde{e}_{i}]=[\nabla\widetilde{e}_{i},J_{\text{M}%
}\widetilde{e}_{i}]+[\widetilde{e}_{i},J_{\text{M}}(\nabla\widetilde{e}_{i})].
\label{1}%
\end{equation}
The definition of $\widetilde{e}_{i}$ implies that $\nabla\widetilde{e}_{i}=0$
for $i=1,...,n.$ Using this fact we obtain:
\begin{equation}
\nabla\lbrack\widetilde{e}_{i},\widetilde{e}_{i+n}]=0. \label{2}%
\end{equation}
Combining formulas \ref{0} and \ref{2} we conclude that
\[
\lbrack\widetilde{e}_{i},\widetilde{e}_{i+n}]\in\mathcal{D}_{i}%
\]
at each point $m\in\mathcal{U}.$ Proposition \ref{g32} is proved.
$\blacksquare$

\begin{proposition}
\label{g321}There exists a unique totally geodesic two dimensional real
submanifold $Z_{i}$ such that at each point $m\in Z_{i}$ we have:
\[
T_{m,Z_{i}}=\mathcal{D}_{i}\left|  _{m}\right.  .
\]

\end{proposition}

\textbf{Proof: }Since the distributions $\mathcal{D}_{i}$ are integrable, we
can apply the Frobenius theorem and deduce the existence of the real
submanifolds $Z_{i}$ such that at each point $m\in Z_{i}$ we have:
$T_{m,Z_{i}}=\mathcal{D}_{i}|_{m}.$ From the way we define the distributions
$\mathcal{D}_{i},$ it follows that each geodesic $\gamma(t)$ which passes
through the point $m_{0}\in$M and
\[
\frac{d\gamma(t)}{dt}\left|  _{t=0}\right.  \in\mathcal{D}_{i}(0)
\]
will be contained in $Z_{i}.$ This is exactly the definition of a totally
geodesic submanifold. $\blacksquare$

Proposition \ref{g321} implies Theorem \ref{g2}. $\blacksquare$

\begin{remark}
\label{g4} Theorem \ref{g2} implies that to the pairs ($e_{i},J_{\text{M}%
}e_{i})$ of vectors in $T_{m_{0},M_{0}}$ we can associate a system of local
holomorphic coordinates ($z^{1},...,z^{n})$ in a small neighborhood of
$m_{0}\in$M$_{0}.$ The complex curves $Z_{i}$ are defined by the system of
equations:
\[
z^{1}=...=z^{i-1}=z^{i+1}=...=z^{n}=0.
\]

\end{remark}

\subsection{The Analogue of Cartan's Formula}

\begin{theorem}
\label{g5} Let g be a K\"{a}hler metric on M$.$ Let ($z^{1},...,z^{n})$ be the
coordinate system as defined in Remark \ref{g4} in a small neighborhood of
$m_{0}\in$M. Then we have the following expression for $g_{i,\overline{j}} $:
\[
g_{i,\overline{j}}(z,\overline{z})=\delta_{i,\overline{j}}+\frac{1}%
{6}R_{i,\overline{j},k,\overline{l}}z^{k}\overline{z}^{l}+...
\]
\textit{where }$R_{i,\overline{j},k,\overline{l}}$ \textit{is the curvature
tensor.}
\end{theorem}

\textbf{Proof:} Let us consider the geodesic coordinates $(x^{1}%
,y^{1};...,x^{n},y^{n})$ associated to the orthonormal frame $(e_{1}%
,J_{\text{M}}e_{1};,,;e_{n},J_{\text{M}}e_{n}).$ In the geodesic coordinate
system associated to $(e_{1},J_{\text{M}}e_{1};,...;e_{n},J_{\text{M}}e_{n}),$
we have:
\[
\frac{\partial}{\partial x^{m}}(g_{i,j})\left|  _{m_{0}}\right.
=\frac{\partial}{\partial y^{m}}(g_{i,j})\left|  _{m_{0}}\right.  =0
\]
for $m=1,...,n.$ This follows from the following Proposition:

\begin{proposition}
\label{g51}At the origin of a geodesic coordinate system, the metric has the
following Taylor expansion:
\[
g_{ij}=\delta_{ij}+\frac{1}{3}R_{ip,gj}x^{p}x^{q}+O(|x|^{3}),
\]
\textit{where }$R_{i,p,q,j}$ \textit{is the curvature tensor of g. (See
}Proposition 1.14 on page 8 of the book \cite{Roe}.)
\end{proposition}

On the other hand, Theorem \ref{g2} implies that we have the following
relations between the coordinate systems $(x^{1},y^{2};...;x^{n},y^{n})$ and
$(z^{1},...,z^{n})$
\[
z^{i}=x^{i}+\sqrt{-1}y^{i}.
\]
So we deduce that
\[
\frac{\partial}{\partial z^{i}}=\frac{1}{2}\left(  \frac{\partial}{\partial
x^{i}}-\sqrt{-1}\frac{\partial}{\partial y^{i}}\right)
\]
and
\[
\frac{\overline{\partial}}{\overline{\partial z^{i}}}=\frac{1}{2}\left(
\frac{\partial}{\partial x^{i}}+\sqrt{-1}\frac{\partial}{\partial y^{i}%
}\right)
\]
where
\[
\frac{\partial}{\partial zk}\left(  g_{i,\overline{j}}\right)  \left\vert
_{m_{0}}\right.  =\frac{\overline{\partial}}{\overline{\partial z^{i}}}\left(
g_{i,\overline{j}}\right)  \left\vert _{m_{0}}\right.  =0.
\]
From here and Proposition \ref{g51} we can conclude that
\[
g_{i,\overline{j}}(z,\overline{z})=\delta_{i,\overline{j}}+\frac{1}%
{6}R_{i,\overline{j},k,\overline{l}}z^{k}\overline{z}^{l}+...
\]
where $R_{i,\overline{j},k,\overline{l}}$ is the curvature tensor. This proves
Theorem \ref{g5}. $\blacksquare$

\begin{corollary}
\label{lines}Let M be a K\"{a}hler Manifold with a K\"{a}hler metric
$g_{i,\overline{j}}(z,\overline{z})$ such that $g_{i,\overline{j}}%
(z,\overline{z})$ are real analytic functions with respect to a local
holomorphic coordinates, then we can identify the totally geodesic one complex
dimensional manifolds $\gamma(\tau)$ in M at a point $m\in$M with respect to
the K\"{a}hler metric $g_{i,\overline{j}}(z,\overline{z})$ with the complex
lines
\[
\mathit{\ }l\subset T_{m,\text{M}}=\mathbb{C}^{N}%
\]
\ through the origin\textit{\ }$0\in\mathbb{C}^{N}$ with directions defined
by\textit{\ }$\frac{d}{d\tau}\gamma(\tau)|_{0}.$\textit{\ }
\end{corollary}

\textbf{Proof: }Let $(z^{1},...,z^{N})$ be\ the holomorphic coordinate system
on $\mathcal{U}$ introduced in Definition \ref{g4}. Theorem \ref{g5} implies
that if we restrict the metric on the line $\tau v_{0}=l,$ we will get the
following formula:
\begin{equation}
(g_{i,\overline{j}}(\tau,\overline{\tau})|_{l}=\delta_{i,\overline{j}}%
+\frac{1}{6}R_{v_{0}}|\tau|^{2}+\text{ higher order terms}, \label{c1}%
\end{equation}
where $R_{v_{0}}$ is the holomorphic sectional curvature in the direction
\[
v_{0}\in T_{m,\text{M}}=\mathbb{C}^{N}.
\]
Formula \ref{c1}, Theorem \ref{g5} and the fact that the K\"{a}hler metric
$g_{i,\overline{j}}(z,\overline{z})$ depends real analytically on the
holomorphic coordinate system $(z^{1},...,z^{N})$ on $\mathcal{U}$ imply
directly Corollary \ref{lines}. $\blacksquare$

\begin{definition}
\label{comexp}The map defined in Cor. \ref{lines} we will be called the
complex exponential map in K\"{a}hler Geometry.
\end{definition}

\section{Appendix II. Applications of Bishop's Convergence Theorem}

We will prove the following Theorem following the arguments in \cite{Siu1}:

\begin{theorem}
\label{BR}Let
\[
\pi:\mathcal{N\rightarrow K}%
\]
and
\[
\pi^{\prime}:\mathcal{N}^{^{\prime}}\rightarrow\mathcal{K}%
\]
be two holomorphic families of K\"{a}hler manifolds over a complex manifold
$\mathcal{K}$ as a parameter space. Suppose that both families can be
identified with the $C^{\infty}$ trivial family $\mathcal{K\times}N_{0}$ by
some diffeomorphism $\phi$ isotopic to the identity$.$ Let $0\in\mathcal{K}$
be a fixed point$.$ Let $\tau_{n}\in\mathcal{K}$ be a sequence of points in
$\mathcal{K}$ such that
\[
\underset{k\rightarrow\infty}{\lim}\tau_{k}=\tau_{0}.
\]
Suppose that for all $k$ there exist biholomorphic maps
\[
\phi_{k}:N_{\tau_{k}}\mathit{\rightarrow}N_{\tau_{k}}^{^{\prime}}%
\]
such that $\phi_{k}$ \textit{induces the identity map on }$H^{2}%
(N_{0},\mathbb{Z}).$ \textit{Then }$N_{\tau_{0}}$\textit{\ and }$N_{\tau_{0}%
}^{^{\prime}}$\textit{\ are bimeromorphic}
\end{theorem}

\textbf{Proof:} For each $\tau_{k}$ we have a biholomorphic map
\[
\phi_{k}:N_{\tau_{k}}^{^{\prime}}\rightarrow N_{\tau_{k}}%
\]
then the graphs of $\phi_{k}$ are submanifolds
\[
\Gamma_{k}\subset N_{\tau_{k}}^{^{\prime}}\times N_{\tau_{k}}.
\]
The idea of the proof is to show that
\begin{equation}
\underset{k\rightarrow\infty}{\lim}\Gamma_{k}=\Gamma_{0}\subset N_{\tau_{k}%
}^{^{\prime}}\times N_{\tau_{0}} \label{Gr}%
\end{equation}
exists and from this fact to deduce Theorem \ref{BR}.

\begin{lemma}
\label{BR11}There exists a subsequence $\{k_{n}\}$ of the sequence $\{k\}$
such that the limit of the currents $[\Gamma_{k_{n}}]$
\[
\underset{n\rightarrow\infty}{\lim}[\Gamma_{k_{n}}]=[\Gamma_{0}]
\]
exists, where $\Gamma_{0}$ is a complex analytic subspace in $N_{\tau_{o}%
}\times N_{\tau_{0}}.$
\end{lemma}

\textbf{Proof:} The proof of Lemma \ref{BR11} is based on the following
Theorem of Bishop. See \cite{Bi}:

\begin{theorem}
\label{Bish}Let M be a complex analytic manifold and let $[N_{n}]$ be a
sequence of complex analytic subspaces considered as currents in M$.$ Suppose
that $h$ is a Hermitian metric on M such that
\[
\text{\textit{vol(}}N_{n})\leq C
\]
\textit{then there exists a subsequence of current }$[N_{n_{k}}]$ such that
the limits of the currents
\[
\underset{k\rightarrow\infty}{\lim}[N_{n_{k}}]=[N_{0}]
\]
\textit{exists and as a current }$[N_{0}]$\textit{\ is defined by a complex
analytic subspace }$N_{0}$\textit{\ in }$M.$
\end{theorem}

In order to apply Theorem \ref{Bish} we need to construct a Hermitian metric
$h$ on the complex manifold
\[
\mathcal{M}=\mathcal{N}^{^{\prime}}\underset{\mathcal{K}}{\times}\mathcal{N}\
\]
and prove that the volume of the graphs $\Gamma_{k}$ of the isomorphisms
$\phi_{k}$ are uniformly bounded .

\textbf{Construction of an Hermitian metric on }$\mathcal{M}.$

We assumed that $\mathcal{N}\rightarrow\mathcal{K}^{^{\prime}}$ and
$\mathcal{N}\rightarrow\mathcal{K}$ are families of K\"{a}hler manifolds. So
we may we have two families of K\"{a}hler metrics $\omega^{1,1}(\tau)$ and
$\omega_{1}^{1,1}(\tau)$ respectively defined on the fibres of $\mathcal{N}%
^{^{\prime}}$ and $\mathcal{N}$ and the metrics depend on $\tau\in\mathcal{K}$
in a $C^{\infty}$ manner. Let
\[
\eta=\sqrt{-1}\sum_{i=1}^{N}d\tau^{i}\wedge\overline{d\tau^{i}}%
\]
be a positive definite form on $\mathcal{U}.$ The collection of forms
$\omega^{1,1}(\tau),$ $\omega_{1}^{1,1}(\tau)$ and $\eta$ define a Hermitian
metric $H$ on $\mathcal{M}=\mathcal{N}^{^{\prime}}\underset{\mathcal{K}%
}{\times}\mathcal{N}$ $.$

\begin{proposition}
\label{BR111} The submanifolds $\Gamma_{k}$ defined by $\left(  \ref{Gr}%
\right)  $ have a constant volume with respect to the family of metrics
$\omega^{1,1}(\tau)+\omega_{1}^{1,1}(\tau)$ on $N_{\tau}\times N_{1,\tau}$,
i.e.
\[
vol(\Gamma_{k})=c
\]
for all $k.$
\end{proposition}

\textbf{Proof:} It is easy to see that
\[
vol(\Gamma_{k})=\int\limits_{N_{k}^{\prime}}(\phi_{k}^{\ast}(\omega_{1}%
^{1,1}(\tau_{k})+\omega^{1,1}(\tau_{k}))^{2n}.
\]
Remember that $\phi_{k}^{\ast}$ is the identity map on $H^{2}(N,\mathbb{Z}).$
Here we are using the fact that as $C^{\infty}$ manifolds $\mathcal{N=K\times
}N_{0}$ and $\mathcal{N}^{\prime}=\mathcal{K}^{^{\prime}}\times N_{0}$ are
diffeomorphic and the diffeomorphism which acts fibrewise is isotopic to the
identity on $N_{0}.$ Using this fact we see that the classes of cohomology of
$\phi_{k}^{\ast}(\omega_{1}^{1,1}(\tau_{k}))$ and $\omega^{1,1}(\tau_{k})$ are
fixed. Let us denote them by $[\omega_{1}]$ and $[\omega]$ respectively. So we
get that
\[
\int\limits_{N_{k}^{\prime}}(\phi_{k}^{\ast}(\omega_{1}^{1,1}(\tau
_{k}))+\omega^{1,1}(\tau_{k}))^{2n}=\int\limits_{N_{k}^{\prime}}([\omega
_{1}]+[\omega])^{2n}=c_{0}.
\]
Proposition \ref{BR111} is proved. $\blacksquare$

\textbf{The end of the proof of Lemma }\ref{BR11}: For a subvariety $Z$ of
pure codimension m in a complex manifold $N_{\tau_{0}},$ we denote by $[Z]$
the current on $N_{\tau_{0}}$ defined by $Z$. Now we can apply Bishop's
Theorem \ref{Bish} and conclude that the sequence of currents converges weakly
to a current $[\Gamma_{0}]$ in $N_{\tau_{0}}\times N_{\tau_{0}}^{^{\prime}}$
of the form
\[
\lbrack\Gamma_{0}]=\sum_{i}m_{i}[\Gamma_{i}],
\]
where $m_{i}$ are positive integers and $\Gamma_{i}$ are irreducible complex
analytic subspaces in $N_{\tau_{0}}\times N_{\tau_{0}}^{^{\prime}}.$ Lemma
\ref{BR11} is proved $\blacksquare$.

\textbf{The end of the proof of Theorem }\ref{BR}: Any closed current $[Z]$ on
$N_{\tau_{0}}\times N_{\tau_{0}}^{^{\prime}}$ defines a linear map:
\[
\lbrack Z]_{\ast}:H^{\ast}(N_{\tau_{0}}^{^{\prime}},\mathbb{C})\rightarrow
H^{\ast}(N_{\tau_{0}},\mathbb{C})
\]
between the cohomology rings of $N_{\tau_{0}}^{^{\prime}}$ and $N_{\tau_{0}}$
as follows: Let
\[
\alpha\in H^{\ast}(N_{\tau_{0}}^{^{\prime}},\mathbb{C}),
\]
then
\[
\lbrack Z]_{\ast}(\alpha):=(pr_{2})_{\ast}\left(  [Z]\wedge(pr_{1})^{\ast
}\alpha\right)
\]
where $pr_{i}$ are respectively the projections of $N_{\tau_{0}}\times
N_{0}^{\prime}$ onto the first and the second \ factor and $\left(
pr_{1}\right)  _{\ast}$ and $\left(  pr_{2}\right)  ^{\ast}$ are the
pushforward and the pullback maps. Similarly we define a linear map
\[
\lbrack Z]^{\ast}:H^{\ast}(N_{\tau_{0}},\mathbb{C})\rightarrow H^{\ast
}(N_{\tau_{0}}^{^{\prime}},\mathbb{C}).
\]
The map $[\Gamma_{0}]_{\ast}$ defined by the current $[\Gamma_{0}]$ clearly
agrees with the map
\[
\underset{k\rightarrow\infty}{\lim}[\Gamma_{k}]=\underset{k\rightarrow\infty
}{\lim}id=id=[\Gamma_{0}]_{\ast}%
\]
on $H^{2n}(N,\mathbb{C}).$ From here we deduce that $[\Gamma_{0}]_{\ast}$ is
the identity map on $H^{2n}(N,\mathbb{C}).$ This implies
\[
\lbrack\Gamma_{0}]_{\ast}(\wedge^{n}[\omega_{1}^{1,1}])=(pr_{2})_{\ast}\left(
\sum_{i}m_{i}[\Gamma_{i}]\wedge(pr_{1})^{\ast}(\wedge^{n}[\omega_{1}%
^{1,1}])\right)  =\wedge^{n}[\omega_{1}^{1,1}].
\]
Hence there must be some $\Gamma_{j}$ among the components of which is
projected both onto $N_{0}^{\prime}$ and $N_{\tau_{0}}.$ There can be only one
$\Gamma_{j}$ in $\Gamma_{0}$ with multiplicity
\[
m_{j}=1,
\]
because both
\[
\left(  \sum_{i}m_{i}[\Gamma_{i}]\right)  _{\ast}and\text{ }\left(  \sum
_{i}m_{i}[\Gamma_{i}]\right)  ^{\ast}%
\]
must leave fixed the class in $H^{0}(N,\mathbb{C})$ defined by the function
$1.$ This implies that the projection maps from $\Gamma_{j}$ to $N_{\tau_{0}}$
and to $N_{\tau_{0}}^{^{\prime}}$ have degree one. So $N_{\tau_{0}}$ and
$N_{\tau_{0}}^{^{\prime}}$ are bimeromorphically equivalent. \textbf{Theorem
}\ref{BR} is proved. $\blacksquare$

\begin{theorem}
\label{BR2}Suppose that the map $f:N_{\tau_{0}}\rightarrow N_{\tau_{0}%
}^{^{\prime}}$ is bimeromorphic. Let $\mathcal{L}_{1}$ be an ample line bundle
on $N_{\tau_{0}}^{^{\prime}}$ such that $f^{\ast}(\mathcal{L}_{1})$ is also an
ample line bundle on $N_{\tau_{0}}.$ Then $f$ is a biholomorphic map.
\end{theorem}

\textbf{Proof:} Since
\[
f:N_{\tau_{0}}\rightarrow N_{\tau_{0}}^{^{\prime}}%
\]
is a bimeromorphic map, there exist complex analytic subspaces $Z_{1}\subset
N_{\tau_{0}}$ and $Z_{2}\subset N_{\tau_{0}}^{^{\prime}}$ of codim $\geq2$
such that
\[
f:N_{\tau_{0}}\circleddash Z_{1}\rightarrow N_{\tau_{0}}^{^{\prime}%
}\circleddash Z_{2}%
\]
is a biholomorphic map between those two Zariski open sets. From here it
follows that $f^{\ast}(\mathcal{L}_{1})$ is a well defined line bundle on
\[
U:=N_{\tau_{0}}\circleddash Z_{1}.
\]
Since $\mathcal{L}_{1}$ is an ample bundle on $N_{\tau_{0}}^{^{\prime}},$
$\left(  \mathcal{L}_{1}\right)  ^{\otimes m}$ will be a very ample bundle on
$N_{\tau_{0}}^{^{\prime}}$ for some positive integer $m.$ Let
\[
\sigma_{0},...,\sigma_{k}\in H^{0}(N_{\tau_{0}}^{^{\prime}},\left(
\mathcal{L}_{1}\right)  ^{\otimes m})
\]
be a basis in $H^{0}(N_{\tau_{0}}^{^{\prime}},\left(  \mathcal{L}_{1}\right)
^{\otimes m})$. Since $\left(  \mathcal{L}_{1}\right)  ^{\otimes m}$ is a very
ample line bundle then $\sigma_{0},...,\sigma_{k}$ define an embedding of
$N_{\tau_{0}}^{^{\prime}}$ into $\mathbb{P}^{k}.$ Clearly $f^{\ast}(\sigma
_{0}),..,f^{\ast}(\sigma_{k})$ will define an embedding of the Zariski open
set
\[
U_{0}:=N_{\tau_{0}}\circleddash Z_{1}%
\]
in $N_{\tau_{0}}$ into $\mathbb{P}^{k}.$ Since $Z_{1}$ has a complex
codimension $\geq2$ Hartogs Theorem implies that the sections
\[
f^{\ast}(\sigma_{0}),..,f^{\ast}(\sigma_{k})
\]
of the line bundle $f^{\ast}(\mathcal{L}_{1})^{\otimes m}$ are well defined
holomorphic sections on $N_{\tau_{0}}.$ We assumed that $f^{\ast}%
(\mathcal{L}_{1})$ is an ample line bundle on $N_{\tau_{0}}$ so we can choose
$m$ to be large enouph to conclude that $f^{\ast}(\mathcal{L}_{1})^{\otimes
m}$ is a very ample line bundle on $N_{\tau_{0}}.$ From here we conclude that
\[
x\in N_{\tau_{0}}\rightarrow(f^{\ast}(\sigma_{0})(x),...,f^{\ast}(\sigma
_{k})(x))
\]
defines a holomorphic embedding
\[
\varphi:N_{\tau_{0}}\subset\mathbb{P}^{k}.
\]
So we deduce that
\[
f:N_{\tau_{0}}\rightarrow N_{\tau_{0}}^{^{\prime}}%
\]
is an isomorphism. Theorem \ref{BR2} is proved. $\blacksquare$

\begin{remark}
The arguments that we used to prove Theorem \ref{BR} are similar to the
arguments used by Siu in \cite{Siu1}. These arguments were suggested by
Deligne to D. Burns and Rapoport.
\end{remark}

\section{Appendix III. Counter Examples to the Analogue of Shafarevich
Conjecture of CY Manifolds}

Let us consider K3 surfaces on which there exists an automorphism $\sigma$ of
order two such that $\sigma$ has no fixed points. Such algebraic K3 surfaces
exist and the quotient $X/\sigma=Y$ is an Enriques surface. It is a well known
fact that $\sigma$ acts on the holomorphic two form $\omega_{X}$ on the K3
surface $X$ as follows:
\begin{equation}
\sigma^{\ast}(\omega_{X})=-\omega_{X}. \label{Enr0}%
\end{equation}
Let us consider an elliptic curve
\[
E_{\tau}:=\mathbb{C}/\Lambda_{\tau},
\]
where
\[
\Lambda_{\tau}:=\left\{  m+n\tau|m,n\in\mathbb{Z},\tau\in\mathbb{C}\text{ and
}\operatorname{Im}\tau>0\right\}  .
\]
We know that $E_{\tau}$ can be embedded in $\mathbb{CP}^{2}$ and in one of the
standard affine open sets $E_{\tau}$ is given by the equation:
\begin{equation}
y^{2}=4x^{3}+g_{2}x+g_{3}, \label{Enr1}%
\end{equation}
where
\[
g_{2}(\tau)=60E_{4}(\tau)\text{ and }g_{3}(\tau)=140E_{6}(\tau).
\]
$E_{4}(\tau)$ and $E_{6}(\tau)$ are the Eisenstein series, i.e. for $n>2$
$E_{2n}(\tau)$ are defined as follows:
\[
E_{2n}(\tau)=\sum_{(m,n)\neq(0,0)}\frac{1}{\left(  m+n\tau\right)  ^{2n}}.
\]
Let $\sigma_{1}$ be the involution defined on the equation $\left(
\ref{Enr1}\right)  $as follows:
\[
\sigma_{1}(x,y)=(-x,y).
\]
Clearly $\sigma_{1}$ acts on $E_{\tau}\,$and
\begin{equation}
E_{\tau}/\sigma_{1}=\mathbb{CP}^{1}. \label{Enr2}%
\end{equation}
From $\left(  \ref{Enr2}\right)  $ we obtain that $\sigma_{1}$ acts on the
holomorphic one form
\[
dz_{\tau}=\frac{dx}{y}%
\]
on $E_{\tau}$ as follows:
\begin{equation}
\sigma_{1}^{\ast}(dz_{\tau})=-dz_{\tau}. \label{Enr3}%
\end{equation}
It is very easy to see that $\left(  \ref{Enr0}\right)  $ and $\left(
\ref{Enr3}\right)  $ imply that the quotient of the product
\[
X\times E_{\tau}/\sigma\times\sigma_{1}=\text{M}%
\]
will be a CY manifold. These CY manifolds are called Borcea Voisin manifolds.
From the theory of moduli of Enriques surfaces and of elliptic curves it
follows that the moduli space $\mathcal{M}$(M) of M is isomorphic to
\[
\Gamma\backslash\mathfrak{h}_{2,10}\times\mathbb{PSL}_{2}(\mathbb{Z}%
\mathbb{)}\backslash\mathfrak{h.}%
\]
For more detail about moduli of Enriques surfaces see \cite{BPV}. Let
$\mathcal{X}\rightarrow C$ be a non isotrivial family of K3 surfaces with an
involution $\sigma$ without fixed points acting on the non singular fibres of
the above family. In \cite{JT98} it was proved that the set $\mathcal{D}_{C}$
of points in $C$ over which the fibres are singular is not empty. We will call
the set $\mathcal{D}_{C}$ the discriminant locus.

Let $E_{\lambda}$ be the family of elliptic curves:
\begin{equation}
y^{2}=x(x-1)(x-\lambda). \label{Enr4}%
\end{equation}
Clearly $\left(  \ref{Enr4}\right)  $ defines a family of elliptic curves over
the projective line. Let us denote it by
\[
\tilde{E}\rightarrow\mathbb{CP}^{1}.
\]
Let us take the product
\begin{equation}
\mathcal{X}\times\tilde{E}\rightarrow C\times\mathbb{CP}^{1}. \label{Enr5}%
\end{equation}
On the family $\left(  \ref{Enr5}\right)  $ we can define the action of
$\sigma$ and $\sigma_{1}.$ By taking the quotient we will get a family
\begin{equation}
\mathcal{Y}\rightarrow C\times\mathbb{CP}^{1}. \label{Enr6}%
\end{equation}
of three dimensional CY manifolds over the product space $C\times
\mathbb{CP}^{1}.$ Clearly that for fixed $\lambda\in\mathbb{CP}^{1}$ and
$\lambda\neq0,1$ \& $\infty$ we can identify all curves
\[
C_{\lambda}=C\times\lambda.
\]
Moreover under this identification the discriminant locuses $\mathcal{D}_{C}%
{}_{\lambda}$ are identified too. This means that the discriminant locus of
the family $\left(  \ref{Enr6}\right)  $ satisfies:
\[
\mathcal{D}_{\mathbb{CP}^{1}\times C}=p_{1}^{\ast}\mathcal{D}_{\mathbb{CP}%
^{1}}+p_{2}^{\ast}\mathcal{D}_{C}.
\]
So the family $\left(  \ref{Enr6}\right)  $ is not rigid and gives a counter
example to Shafarevich's conjecture.

\end{document}